\magnification=\magstep1
\input amstex
\documentstyle{amsppt}

\define\defeq{\overset{\text{def}}\to=}
\define\ab{\operatorname{ab}}
\define\pr{\operatorname{pr}}
\define\Gal{\operatorname{Gal}}

\def \isom {\overset \sim \to \rightarrow}

\define\Pic{\operatorname{Pic}}

\define\Spec{\operatorname{Spec}}
\define\id{\operatorname{id}}

\define\Ker{\operatorname{Ker}}

\def \c{\operatorname {c}}

\def \et{\operatorname {et}}

\def\char{\operatorname{char}}

\def\Sect{\operatorname{Sect}}
\def\genus{\operatorname{genus}}
\def\sep{\operatorname{sep}}
\def\and{\operatorname{and}}
\def\sol{\operatorname{sol}}
\def\Br{\operatorname{Br}}
\def\cl{\operatorname{cl}}

\define\Primes{\frak{Primes}}
\NoRunningHeads
\NoBlackBoxes
\topmatter

\title
The Cuspidalisation of Sections of Arithmetic Fundamental Groups II
\endtitle

\author
Mohamed Sa\"\i di
\endauthor

\abstract In this paper, which is a sequel to [Sa\"\i di], we investigate the theory of cuspidalisation of sections of arithmetic fundamental groups of hyperbolic curves
to cuspidally $i$-th and $2/p$-th step prosolvable 
arithmetic fundamental groups. As a consequence we exhibit two, necessary and sufficient, conditions 
for sections of arithmetic fundamental groups of hyperbolic curves over $p$-adic local fields to arise from rational points. We also exhibit a class of sections of arithmetic 
fundamental groups of $p$-adic curves which are orthogonal to $\Pic^{\wedge}$, and which satisfy (unconditionally) one of the above conditions.
\endabstract
\toc

\subhead
\S 0. Introduction
\endsubhead

\subhead
\S 1. Cuspidally $i$-th and $i/p$-th step prosolvable geometric fundamental groups
\endsubhead

\subhead
\S 2. Cuspidalisation of sections of cuspidally $i$-th step prosolvable arithmetric fundamental groups
\endsubhead

\subhead
\S 3. Lifting of sections to cuspidally $2/p$-th step prosolvable arithmetric fundamental groups
\endsubhead

\subhead
\S 4. Geometric sections of arithmetic fundamental groups of $p$-adic curves
\endsubhead

\subhead
\S 5. Local sections of arithmetic fundamental groups of $p$-adic curves 
\endsubhead


\endtoc

\endtopmatter

\document

\subhead
\S 0. Introduction
\endsubhead
Let $k$ be a characteristic zero field, $X$ a proper, smooth, and geometrically connected hyperbolic (i.e., $\genus (X)\ge 2$) algebraic curve over
$k$. Let $K_X$ be the function field of $X$, $K_X^{\sep}$ a separable closure of $K_X$, and $\overline k\subset K_X^{\sep}$ the algebraic closure of $k$.
Let $\pi_1(X)$ be the \'etale fundamental group of $X$ which sits in the following exact sequence
$$1\to \pi_1(\overline X)\to \pi_1(X) @>{\pr}>> G_k\to 1,$$
where $G_k\defeq \Gal (\overline k/k)$, and $\pi_1(\overline X)$ is the geometric \'etale fundamental group of $X$.
Let $G_X\defeq \Gal (K_X^{\sep}/K_X)$, and $\overline G_X\defeq \Gal (K_X^{\sep}/K_X.\overline k)$.
Thus, we have exact sequences
$$1\to \overline G_X\to G_X\to G_k\to 1,$$
and
$$1\to \Cal I_X\to G_X\to \pi_1(X)\to 1,$$
where $\Cal I_X$ is the inertia subgroup.
The theory of cuspidalisation of sections of arithmetic fundamental groups was initiated in [Sa\"\i di], its ultimate aim is to reduce the Grothendieck anabelian 
section conjecture to its birational version.  It can be formulated as follows (cf. loc. cit.).
\definition
{The Cuspidalisation Problem for Sections of $\pi_1(X)$}
Let $G_X\twoheadrightarrow H\twoheadrightarrow \pi_1(X)$ be a quotient of $G_X$.
Given a section $s:G_k\to \pi_1(X)$
of the projection $\pi_1(X)\twoheadrightarrow G_k$, is it possible to {\bf lift} $s$ to a section  $\tilde s:G_k\to H$
of the projection $H\twoheadrightarrow G_k$? i.e., is it possible to construct a section $\tilde s$ such that
the following diagram is commutative
$$
\CD
G_k @>{\tilde s}>> H \\
@V{\id}VV   @VVV \\
G_k @>{s}>>  \pi _1(X)
\endCD
$$
where the right vertical map is the projection $H\twoheadrightarrow \pi_1(X)$?
\enddefinition

In [Sa\"\i di] we investigated the cuspidalisation problem in the case $H\defeq G_X^{(\c-\ab)}$ is the maximal (geometrically) cuspidally abelian quotient of $G_X$. 
In this paper we generalise this theory to the (geometrically) {\it cuspidally $i$-th} as well as $i/p$-th; where $p$ is a prime, {\it step prosolvable} quotient of $G_X$.

For $i\ge 0$, let $\Cal I_X\twoheadrightarrow \Cal I_{X,i}$ be the maximal i-th step prosolvable quotient of $\Cal I_X$, and
$G_{X}^{(i-\sol)}\defeq G_X/\Ker (\Cal I_X \twoheadrightarrow \Cal I_{X,i})$ the maximal (geometrically) cuspidally i-th step prosolvable quotient of $G_X$ 
($G_X^{(\c-\ab)}\defeq G_{X}^{(1-\sol)}$). 
For $i\ge 1$, let $s_i:G_k\to  G_{X}^{(i-\sol)}$ be a section of the projection $G_{X}^{(i-\sol)}\twoheadrightarrow G_k$. In $\S2$ we investigate the problem of lifting $s_i$ to a section
 $s_{i+1}:G_k\to  G_{X}^{(i+1-\sol)}$ of the projection $G_{X}^{(i+1-\sol)}\twoheadrightarrow G_k$. 
We say that the field $k$ satisfies the condition $\bold {(H)}$ if the following holds. The Galois cohomology groups $H^1(G_k,M)$
are {\it finite} for every finite $G_k$-module $M$. This condition is satisfied for instance if the Galois group $G_k$ is (topologically) finitely generated (e.g. $k$ is a $p$-adic local field).
One of our main results in this paper is the following (cf. Theorem 2.3.8).

\proclaim{Theorem 1} Assume ${i\ge 1}$,
and $k$ satisfies the condition $\bold {(H)}$. 
The section $s_i:G_k\to G_X^{(i-\sol)}$ {\bf lifts} to a section $s_{i+1} : G_k\to G_X^{(i+1-\sol)}$
of the projection $G_X^{(i+1-\sol)}\twoheadrightarrow G_k$ {\bf if and only if} for every $X'\to X$ 
a neighbourhood of the section $s_i$ (i.e., corresponding to an open subgroup of $G_X^{(i-\sol)}$
containing $s_i(G_k)$) the class of $\Pic _{X'}^1$ in $H^1(G_k,\Pic_{X'}^0)$ 
lies in the maximal divisible subgroup of  $H^1(G_k,\Pic _{X'}^0)$.
\endproclaim
 
 Key to the proof of Theorem 1 is the description of the $G_k$-module structure, induced by $s_i$, of $\Cal I_X[i+1]\defeq \Ker(G_X^{(i+1-\sol)}\twoheadrightarrow G_X^{(i-\sol)})$
as the projective limit of the Tate modules of the jacobians of the neighbourhoods $\{X'\}$ as in the statement of Theorem 1 (cf. Proposition 1.1.5, and Lemma 2.3.2).
 
 In $\S3$ we investigate the following mod-$p$ variant of Theorem 1, where $p$ is a prime integer. 
 Let $t\ge 0$, $i\ge 0$, and $\Cal I_X\twoheadrightarrow \Cal I_{X,i+1/p^t}$ the $i+1$-th
 quotient of the $\Bbb Z/p^t\Bbb Z$-derived series of $\Cal I_X$ (cf. 1.2). Thus, $\Cal I_{X,i+1/p^t}$ is $i+1$-step prosolvable with successive abelian quotients annihilated by $p^t$.
Let $G_X^{(i+1/p^t-\sol)}\defeq G_X/\Ker (\Cal I_X\twoheadrightarrow \Cal I_{X,i+1/p^t})$ be the maximal (geometrically) cuspidally $i+1/p^t$-th prosolvable quotient of $G_X$.
Given a section $s:G_k\to \pi_1(X)$ of the projection $\pi_1(X)\twoheadrightarrow G_k$ we investigate the problem of lifting $s$ to a section
$s_{i+1}:G_k\to G_X^{(i+1/p^t-\sol)}$ of the projection $G_X^{(i+1/p^t-\sol)}\twoheadrightarrow G_k$ in the {\bf case} $\bold {i=1}$ and $\bold {t=1}$; the only case needed for applications in 
$\S4$, and $\S5$. 
For this purpose we introduce a certain quotient $G_X\twoheadrightarrow G_X^{(p,2)}\twoheadrightarrow G_X^{(2/p-\sol)}$, and investigate 
the problem of lifting the section $s$ to a section
$\tilde s:G_k\to G_X^{(p,2)}$ of the projection $G_X^{(p,2)}\twoheadrightarrow G_k$ (this would give rise to a section of the projection $G_X^{(2/p-\sol)}\twoheadrightarrow G_k$
which lifts $s$). 
The quotient $G_X^{(p,2)}$ sits in an exact sequence $1\to \Cal I_X[p,2]\to G_X^{(p,2)}\to G_X^{(1/p^2-\sol)}\to 1$, where $\Cal I_X[p,2]$ is abelian annihilated by $p$ 
(cf. 3.3 for more details). In Theorem  3.4.11 we give necessary and sufficient conditions for the section $s$ to lift to a section $\tilde s:G_k\to G_X^{(p,2)}$ (cf. loc. cit. for 
a more precise statement).

In $\S4$ we assume $k$ is a $p$-adic local field (finite extension of $\Bbb Q_p$). We observe in this case that if $\tilde s:G_k\to G_X^{(p,2)}$ is a section of the projection $G_X^{(p,2)}\twoheadrightarrow G_k$, and $s:G_k\to \pi_1(X)$ is the induced section of $\pi_1(X)\twoheadrightarrow G_k$, then $s$ is {\it geometric} in the sense that it arises 
from a rational point $x\in X(k)$ (cf. Proposition 4.6).
Further, we provide the following characterisation of sections $s:G_k\to \pi_1(X)$ which are geometric (cf. Theorem 4.5 where we prove a pro-$\Sigma$; $p\in \Sigma$ is a set of primes, variant of Theorem 2).

\proclaim {Theorem 2} Assume $k$ is a $p$-{\bf adic local field}. A section $s:G_k\to \pi_1(X)$  of the projection $\pi_1(X)\twoheadrightarrow G_k$
 is {\bf geometric} (cf. Definition 4.1) if and only if the following two conditions hold.

\noindent
(i)\ The section $s$ has a {\bf cycle class uniformly orthogonal to $\Pic$ mod-$p^2$} (cf. Definition 3.4.1).
 
 \noindent
(ii)\  There {\bf exists} a section $s':G_k\to G_{X}^{(1/p^2-\sol)}$ 
of the projection  $G_{X}^{(1/p^2-\sol)}  \twoheadrightarrow G_k$ which {\bf lifts} the section $s$ (this holds if condition (i) is satisfied by Theorem 3.4.4) 
such that the following holds. For every
$X'\to X$ a neighbourhood of the section $s'$ (i.e., corresponding to an open subgroup of $G_X^{(1/p^2-\sol)}$containing $s'(G_k)$) the class of $\Pic _{X'}^1$ in $H^1(G_k,\Pic_{X'}^0)$ 
is {\bf divisible by} $p$.
\endproclaim

Condition (ii) in Theorem 2 is a necessary and sufficient condition for the section $s'$ therein to lift to a section
$\tilde s:G_k\to G_X^{(p,2)}$ (cf. Theorem 3.4.10). 

As an application of Theorem 2 we prove the following $p$-adic absolute anabelian result (cf. Theorem 4.8).

\proclaim{Theorem 3} Let $p_X,p_Y$ be prime integers, and
$X$ (resp. $Y$) a proper, smooth, geometrically connected hyperbolic curve over a $p_X$-adic local field
$k_X$ (respectively,  $p_Y$-adic local field $k_Y$). Let $p_X\in \Sigma_X$ (resp. $p_Y\in \Sigma_Y$) be a non-empty set of prime integers of cardinality $\ge 2$, $\Pi_X$ (resp. $\Pi_Y$)
the geometrically pro-$\Sigma_X$ (resp. pro-$\Sigma_Y$) arithmetic fundamental group of $X$ (resp. $Y$),
and $\varphi:\Pi_X\to \Pi_Y$ an isomorphism of profinite groups which fits in the following commutative diagram
$$
\CD
G_X^{(p_X,2)} @>{\widetilde \varphi}>> G_Y^{(p_Y,2)}\\
@VVV   @VVV\\
\Pi_X@>{\varphi}>>\Pi_Y\\
\endCD
$$
where $\widetilde \varphi$ is an isomorphism of profinite groups. Here $G_X^{(p_X,2)}$ (resp. $G_Y^{(p_Y,2)}$) is the pro-$\Sigma_X$ (resp. $\Sigma_Y$) version of the 
profinite group $G_X^{(p_X,2)}$ (resp. $G_Y^{(p_Y,2)}$) (cf. 3.3), and the vertical maps are the natural projections. Then $\varphi$ is geometric, i.e., arises from a uniquely determined isomorphism of schemes $X\isom Y$.
\endproclaim

Finally, in $\S5$ we investigate {\it local} sections of arithmetic fundamental groups of $p$-adic curves. These are sections which arise from sections of 
arithmetic fundamental groups of {\it formal fibres} (cf. Definition 5.2). A geometric section is necessarily local in this sense. Our main result is 
the following (cf. Theorem 5.3).

\proclaim{Theorem 4} Assume $k$ is a {\bf $p$-adic local field}, and $s:G_k\to \pi_1(X)$ is a {\bf local} section
of the projection $\pi_1(X)\twoheadrightarrow G_k$. Then $s$ has {\bf a cycle class which is uniformly orthogonal to $\Pic^{\wedge}$} in the sense of [Sa\"\i di], Definition 1.4.1(i). 
\endproclaim

To the best of our knowledge, local sections of arithmetic fundamental groups of $p$-adic curves are the first non trivial (i.e., not known to be geometric a priori) examples of sections 
of arithmetic fundamental groups of $p$-adic curves which are orthogonal to $\Pic^{\wedge}$. In particular, local sections satisfy condition (i) in Theorem 2. 

\subhead
Notations
\endsubhead
Throughout this paper $\Primes$ denotes the set of all prime integers. 
For a profinite group $H$, we write $H^{\ab}$ for the maximal abelian quotient of $H$.

\subhead
\S1. Cuspidally $i$-th and $i/p$-th step prosolvable geometric  fundamental groups
\endsubhead
Let $\ell $ be an algebraically closed field  of characteristic $l\ge 0$, $X$ a proper smooth and connected
hyperbolic curve (i.e., $\genus(X)\ge 2$)  over $\ell$, and $K_X$ its function field.
Let $\eta $ be a geometric point of $X$ above
its generic point; which determines a separable closure $K_X^{\sep}$ of $K_X$, and
$\pi_1(X,\eta)$ the \'etale fundamental group of $X$ with base
point $\eta$. 

Let 
$\emptyset\neq \Sigma \subseteq \Primes$
be a set of prime integers.
In case $\char (\ell)=l >0$ we assume that $l\notin \Sigma$. 
Write
$\Delta_X\defeq \pi_1(X,\eta)^{\Sigma}$
for the maximal pro-$\Sigma$ quotient of $\pi_1(X, \eta)$. 
Let $\{x_s\}_{s=1}^n\subset X(\ell)$,
$U\defeq X\setminus \{x_1,\cdots,x_n\}$ an open subscheme of $X$,
$\Delta_U\defeq \pi_1(U,\eta)^{\Sigma}$
the maximal pro-$\Sigma$ quotient of the \'etale fundamental group $\pi_1(U,\eta)$ of $U$
with base point $\eta$, and
$I_U\defeq \Ker (\Delta_U\twoheadrightarrow \Delta_X)$.
We shall refer to $I_U$ as the {\it cuspidal subgroup} of $\Delta_U$ with respect to the natural projection $\Delta_U\twoheadrightarrow \Delta _X$
(cf. [Mochizuki], Definition 1.5); it is the
subgroup of $\Delta_U$ (normally) generated by the (pro-$\Sigma$) inertia subgroups at the points 
$\{x_i\}_{i=1}^n$.
We have the following exact sequence
$$1\to I_U\to \Delta_U\to \Delta_X\to 1.$$

\subhead {1.1 The quotient $\Delta_{U,i}$}
\endsubhead
For a profinite group $H$, we denote by $\overline {[H,H]}$ the closed subgroup of $H$ topologically generated by the commutator subgroup.
Consider the derived series of $I_U$
$$.....\subseteq I_U(i+1)\subseteq I_U(i)\subseteq......\subseteq I_U(1)\subseteq I_U(0)=I_U, \tag 1.1$$
where, for $i\ge 0$,
$I_U(i+1)\defeq \overline {[I_U(i),I_U(i)]}$
is the $(i+1)$-th derived subgroup, which is a characteristic subgroup of $I_U$. Write
$$I_{U,i}\defeq I_U/I_U(i).$$
Thus, $I_{U,i}$ is the maximal $i$-th step prosolvable quotient of $I_U$, and $I_{U,1}$ is the maximal abelian quotient of $I_U$. 
There exists a natural exact sequence
$$ 1\to   I_U[i+1]\to I_{U,i+1}\to I_{U,i}\to 1 \tag 1.2$$
where $I_U[i+1]$ is the subgroup $I_U(i)/I_U(i+1)$ of $I_{U,i+1}$ and $I_U[i+1]$ is abelian.
Write 
$$\Delta _{U,i}\defeq \Delta_U/I_{U}(i).$$
We shall refer to $\Delta _{U,i}$ 
(resp. $\Delta _{U,1}$) as the maximal {\bf cuspidally i-th step prosolvable} (resp. maximal {\bf cuspidally abelian})
quotient of $\Delta _U$ (with respect to the surjection $\Delta _U\twoheadrightarrow \Delta _X$).
We have the following commutative diagram of exact sequences.

$$
\CD
 @. 1 @.  1@. \\
@. @VVV    @VVV\\
@. I_U[i+1]  @= I_U[i+1]\\
@.  @VVV   @VVV \\
1 @>>> I_{U,i+1}  @>>>  \Delta _{U,i+1} @>>> \Delta _X @>>> 1\\
@.  @VVV     @VVV   @V{\id_{\Delta_X}}VV\\
1 @>>> I_{U,i}  @>>>  \Delta _{U,i} @>>> \Delta _X @>>> 1\\
@.  @VVV   @VVV \\
@. 1@. 1@.\\
\endCD
\tag 1.3$$

The profinite group $\Delta _{U,i}$; being a quotient of $\Delta_U$, is topologically finitely generated (cf. [Grothendieck], Expos\'e X, Corollaire 3.10, recall $\char (\ell)\notin \Sigma$). 
Hence there exists a sequence of characteristic open subgroups
$$...\subseteq \Delta _{U,i}[j+1]\subseteq \Delta _{U,i}[j]\subseteq...\subseteq \Delta _{U,i}[1]\defeq \Delta _{U,i}$$
of $\Delta_{U,i}$ such that $\bigcap _{j\ge 1}\Delta _{U,i}[j]=\{1\}$. 
The open subgroup $\Delta _{U,i}[j]\subseteq \Delta_{U,i}$ corresponds to a
finite (Galois) cover  $X_{i,j}^U\to X$ between smooth connected and proper $\ell$-curves which is \'etale above $U$. 
The geometric point $\eta$ determines naturally a geometric point $\eta_{i,j}$ of $X_{i,j}^U$. Write
$\Delta_{i,j}^U\defeq\Delta_{X_{i,j}^U}\defeq \pi_1(X_{i,j}^U,\eta_{i,j})^{\Sigma}$
for the maximal pro-$\Sigma$ \'etale fundamental group of $X_{i,j}^U$ with base point $\eta_{i,j}$, and 
$(\Delta_{i,j}^U)^{\ab}$
for the maximal abelian quotient of $\Delta_{i,j}^U$. 
The following Proposition provides a description of the structure of the profinite group $I_U[i+1]$ (cf. sequence (1.2) and diagram (1.3)) in the case
${i\ge 1}$. A description of the structure of $I_U[1]$
is given in [Mochizuki] Proposition 1.14 (see also [Sa\"\i di], 2.1).

\proclaim {Proposition 1.1.1} Let ${i\ge 1}$. There exists a natural isomorphism 
$$I_U[i+1]\isom \underset{j\ge 1} \to{\varprojlim}\ (\Delta_{i,j}^U)^{\ab}.$$
\endproclaim

\demo {Proof of Proposition 1.1.1}
Let $G$ be a finite quotient of $\Delta_{U,i+1}$, which inserts in 
the following commutative diagram
$$
\CD
1 @>>> I_{U,i+1}  @>>>  \Delta _{U,i+1} @>>> \Delta _X @>>> 1\\
@. @VVV     @VVV   @VVV\\
1 @>>> I  @>>> G @>>> G^{\et} @>>> 1\\
\endCD
$$
where the vertical maps are surjective.
We assume, without loss of generality, that $G$ is {\it not} a quotient of $\Delta _{U,i}$.
The quotient $G$ corresponds to a finite Galois cover $X_1\to X$ with Galois group $G$, which factorizes as
$X_1\to X_1^{\et}\to X$, where $X_1^{\et}\to X$ is the maximal \'etale sub-cover with Galois group $G^{\et}$, and $X_1\to X_1^{\et}$ is a (tamely) 
ramified Galois cover with group $I$. 
For $s\in \{1,\cdots,n\}$, let $I_{x_s}\subset G$ be an inertia subgroup associated to $x_s$. Thus, $I_{x_s}$ is only defined up to conjugation, and $I$ is an $(i+1)$-th step 
solvable group (normally) generated by the $I_{x_s}$'s.
Moreover, $I_{x_s}$ is cyclic of order $e_s\ge 1$ (coprime to $l=\char(\ell )$) as the ramification is tame. 
The following claim follows immediately from the well-known structure of $\Delta _U$ (cf. [Grothendieck], Expos\'e X, Corollaire 3.10). 

\proclaim {Claim 1.1.2} There exists a finite quotient $G'$  of $\Delta_{U,i}$, which inserts in the following commutative diagram (where the vertical maps are surjective)
$$
\CD
1 @>>> I_{U,i}  @>>>  \Delta _{U,i} @>>> \Delta _X @>>> 1\\
@. @VVV     @VVV   @VVV\\
1 @>>> I'  @>>> G' @>>> {G '}^{\et} @>>> 1\\
\endCD
$$
such that the following holds.
The quotient $\Delta _{U,i}\twoheadrightarrow G'$ 
corresponds to a finite Galois cover $X_2\to X$ with Galois group $G'$, which factorizes as
$X_2\to X_2^{\et}\to X$, where $X_2^{\et}\to X$ is the maximal \'etale sub-cover with Galois group ${G'}^{\et}$, and $X_2\to X_2^{\et}$ is a (tamely) ramified Galois cover with 
Galois group $I'$ (an $i$-th step solvable group). Further,
for $s\in \{1,\cdots,n\}$, $I'_{x_s}\subseteq I'$ an inertia subgroup associated to $x_s$, 
then $I'_{x_s}$ is {\it cyclic} of order $f_s=e_sh_s$ a multiple of $e_s$.
\endproclaim

Next, let $K_1\defeq K_{X_1}$ (resp. $K_2\defeq K_{X_2}$) be the function field of $X_1$ (resp. $X_2$). 
Let $L\defeq K_1.K_2$ be the compositum of $K_1$ and $K_2$ (in $K_X^{\sep}$), and $\widetilde X$ the normalisation of $X$ in $L$. Thus, $\widetilde X\to X$ is a Galois cover with Galois group $H\subseteq G\times G'$ which is \'etale above $U$ and factorizes as $\widetilde X\to {\widetilde X}^{\et}\to X$, where ${\widetilde X}^{\et}\to X$ is the maximal \'etale sub-cover with Galois group $H^{\et}$, and $\widetilde X\to {\widetilde X}^{\et}$ is a (tamely) ramified Galois cover with group $I_H$
: the subgroup of $H$ (normally) generated by the inertia subgroups at the points of $\widetilde X$ above the $\{x_s\}_{s=1}^n$.
(Thus, we have an exact sequence $1\to I_H\to H\to H^{\et}\to 1$.)

\proclaim {Lemma 1.1.3}  The quotient $\Delta _U\twoheadrightarrow H$ factorizes as   $\Delta _U\twoheadrightarrow \Delta _{U,i+1}\twoheadrightarrow H$. 
\endproclaim

\demo{Proof of Lemma 1.1.3} Indeed, one verifies easily that $I_H$ is a subgroup of $I\times I'$ and $I\times I'$ is $(i+1)$-th step solvable.
\qed
\enddemo

Next, let $\widetilde I$  be the maximal $i$-th step solvable quotient of $I$,
which inserts in the exact sequence
$1\to I(i+1)\to I\to \widetilde I\to 1$, with $I(i+1)$ abelian (note that $I(i+1)$ is non trivial by our assumption that  $G$ is not a quotient of $\Delta_{U,i}$). 
Write $\widetilde G\defeq G/I(i+1)$,
which inserts in the exact sequence $1\to \widetilde I \to  \widetilde G \to G^{\et} \to 1$. In particular,
$\widetilde G$ is a quotient of $\Delta _{U,i}$.
Let $\widetilde H$ be the image of $H$ in $\widetilde G\times G'$. We have a commutative diagram of exact sequences where the vertical maps are natural inclusions.
$$
\CD
1  @>>> I_H(i+1)\defeq H\cap \lgroup I(i+1)\times \{1\} \rgroup @>>> H @>>> \widetilde H @>>>1\\
@. @VVV  @VVV  @VVV  \\
1 @>>> I(i+1)\times \{1\}  @>>> G\times G' @>>> \widetilde G\times G' @>>> 1\\
\endCD
$$

\proclaim {Lemma 1.1.4}
The group $\widetilde H$ is a quotient of $\Delta _{U,i}$. Moreover, the cover $\widetilde X\to X$ factorizes as $\widetilde X\to \widetilde X'\to X$,
where $\widetilde X'\to X$ is Galois (\'etale above $U$) with Galois group $\widetilde H$, and $\widetilde X\to \widetilde X'$ is an {\bf abelian \'etale} cover with Galois group $I_H(i+1)$.
\endproclaim

\demo{Proof of Lemma 1.1.4} The first assertion follows from the various definitions. Next,
the Galois cover $X_1\to X$ factorizes as $X_1\to \widetilde X_1\to X$ where $X_1\to \widetilde X_1$ is Galois with Galois group $I(i+1)$, and 
$\widetilde X_1\to X$ is Galois with group $\widetilde G$. Let $\widetilde X'$ be the normalisation of $X$ in the compositum of the function fields of  $\widetilde X_1$ and $X_2$.
Thus, $\widetilde X'\to X$ is a Galois cover with Galois group $\widetilde H$, and
we have the following commutative diagram
$$
\CD
\widetilde X @>>> X_1\\
@VVV     @VVV\\
\widetilde X' @>>>  \widetilde X_1\\
@VVV   @VVV\\
X_2   @>>> X\\
\endCD
$$
of finite Galois covers.
The ramification index in the Galois cover $X_2 \to X$ 
above a branched (closed) point $x_s\in X$ is divisible by the ramification index above $x_s$ in the Galois cover $X_1\to X$ (cf. the above condition 
that $f_s$ is divisible by $e_s$). The fact that the morphism $\widetilde X\to X_2$, and a fortiori $\widetilde X\to \widetilde X'$; which is abelian with Galois group $I_H(i+1)$, is \'etale follows from Abhyankar's Lemma (cf. [Grothendieck], Expos\'e X, Lemma 3.6).
\qed
\enddemo

Going back to the proof of Proposition 1.1.1, the above discussion shows that the finite quotients $\Delta _{U,i+1}\twoheadrightarrow H$ as in Lemma 1.1.3 form a cofinal system 
of finite quotients of $\Delta _{U,i+1}$. Thus,
$\Delta _{U,i+1}\isom \underset{H} \to{\varprojlim}\ H$. Proposition 1.1.1 then follows from the facts that the various $H$ above fit in an exact sequence
$1\to I_H(i+1)\to H\to \widetilde H\to 1$; $\Delta _{U,i}\isom \underset{\widetilde H} \to{\varprojlim}\ \widetilde H$, and the above Galois covers $\widetilde X\to \widetilde X'$ 
with group $I_H(i+1)$ are \'etale abelian (cf. Lemma 1.1.4).
This finishes the proof of Proposition 1.1.1.
\qed
\enddemo

Similarly, let $G_{K_X}\defeq \Gal (K_X^{\sep}/K_X)$, and $G_X\defeq G_{K_X}^{\Sigma}$ the maximal pro-$\Sigma$
quotient of $G_{K_X}$. We have a natural exact sequence
$$1\to \Cal I_X\to G_X\to \Delta_X\to 1,$$ 
where $\Cal I_X\defeq \Ker (G_X\twoheadrightarrow \Delta_X)$ is the  {\it cuspidal subgroup} of $G_X$ (with respect to the surjection $G_X\twoheadrightarrow \Delta _X$). 
Let $i\ge 0$ and write
$$\Cal I_{X,i}\defeq \Cal I_X/\Cal I_X(i).$$
Thus, $\Cal I_{X,i}$ is the maximal $i$-th step prosolvable quotient of $\Cal I_X$, and $\Cal I_{X,1}$ is the maximal abelian quotient of $\Cal I_X$. 
There exists a natural exact sequence
$$ 1\to   \Cal I_X[i+1]\to \Cal I_{X,i+1}\to \Cal I_{X,i}\to 1, \tag 1.4$$
where $\Cal I_X[i+1]$ is the subgroup $\Cal I_X(i)/\Cal I_X(i+1)$ of $\Cal I_{X,i+1}$, and $\Cal I_X[i+1]$ is abelian.
Write 
$$G_{X,i}\defeq G_X/\Cal I_{X}(i).$$
We shall refer to $G _{X,i}$ 
(resp. $G_{X,1}$) as the maximal {\bf cuspidally i-th step prosolvable} (resp. maximal {\bf cuspidally abelian})
quotient of $G_X$ (with respect to the surjection $G_X\twoheadrightarrow \Delta _X$).
We have the following commutative diagram 

$$
\CD
 @. 1 @.  1@. \\
@. @VVV    @VVV\\
@. \Cal I_X[i+1]  @= \Cal I_X[i+1]\\
@.  @VVV   @VVV \\
1 @>>> \Cal I_{X,i+1}  @>>>  G_{X,i+1} @>>> \Delta _X @>>> 1\\
@.  @VVV     @VVV   @VVV\\
1 @>>> \Cal I_{X,i}  @>>>  G _{X,i} @>>> \Delta _X @>>> 1\\
@.  @VVV   @VVV \\
@. 1@. 1@.\\
\endCD
\tag 1.5$$
of exact sequences.

\proclaim {Proposition 1.1.5} There are natural isomorphisms $G_{X,i}\isom \underset{U} \to{\varprojlim}\ \Delta _{U,i}$, $\Cal I_{X,i}\isom \underset{U} \to{\varprojlim}\ I _{U,i}$,
and $\Cal I_{X}[i+1]\isom \underset{U} \to{\varprojlim}\ I _{U}[i+1]$, where the projective limit is taken over all non-empty open subschemes $U\subseteq X$. Moreover, for ${i\ge 1}$,
we have a natural isomorphism
$$\Cal I_{X}[i+1]\isom \underset{U} \to{\varprojlim}     \lgroup \underset{j\ge 1} \to{\varprojlim}  ({\Delta _{i,j}^{U}})^{\ab} \rgroup,$$
where $({\Delta _{i,j}^{U}})^{\ab}$ is as in the discussion preceding Proposition 1.1.1.
\endproclaim

\demo{Proof} Follows from the various definitions and Proposition 1.1.1.
\qed
\enddemo

\subhead {1.2. The quotient  $\Delta_U\twoheadrightarrow \Delta_U^{p,i+1}$}
\endsubhead
In this subsection we discuss a certain variant of the theory in 1.1, we use the same notations as in loc. cit..
For a profinite group $H$, a prime integer $p$, and an integer $t\ge 1$, write
$$.....\subseteq H(i+1/p^t)\subseteq H(i/p^t)\subseteq......\subseteq H(1/p^t)\subseteq H(0/p^t)=H$$
for the $\Bbb Z/p^t\Bbb Z$-derived series of $H$, where
$H(i+1/p^t)\defeq \overline {<[H(i/p^t),H(i/p^t)], H(i/p^t)^{p^t}>}$
is the $i+1/p^t$-th derived subgroup, which is a characteristic subgroup of $H$. Write
$$H_{i/p^t}\defeq H/H(i/p^t).$$
We will refer to $H_{i/p^t}$ as the maximal $i/p^t$-th step {\it prosolvable} quotient of $H$, and $H_{1/p^t}$ as the maximal {\it abelian} annihilated by $p^t$ quotient of $H$. 
There exists a natural exact sequence
$$ 1\to   H[i+1/p^t]\to H_{i+1/p^t}\to H_{i/p^t}\to 1 $$
where $H[i+1/p^t]$ is the subgroup $H(i/p^t)/H(i+1/p^t)$ of $H_{i+1/p^t}$, and $H[i+1/p^t]$ is abelian annihilated by $p^t$.

Next, let
$p\in \Sigma$,
and consider the $\Bbb Z/p\Bbb Z$-derived series of $I_U$
$$.....\subseteq I_U(i+1/p)\subseteq I_U(i/p)\subseteq......\subseteq I_U(1/p)\subseteq I_U(0/p)=I_U\tag 1.6$$
(cf. the above discussion in the case $t=1$). Then
$I_{U,i/p}\defeq I_U/I_U(i/p)$
is the maximal $i/p$-th step prosolvable quotient of $I_U$, and $I_{U,1/p}$ is the maximal abelian annihilated by $p$ quotient of $I_U$. 
Write 
$\Delta _{U,i/p}\defeq \Delta_U/I_{U}(i/p)$,
which inserts in the exact sequence
$$1\to I_{U,i/p} \to  \Delta _{U,i/p} \to \Delta _X \to 1.$$
We shall refer to $\Delta _{U,i/p}$ 
(resp. $\Delta _{U,1/p}$) as the maximal {\bf cuspidally $i/p$-th step prosolvable} (resp. maximal {\bf cuspidally abelian annihilated by $p$})
quotient of $\Delta _U$ (with respect to the surjection $\Delta _U\twoheadrightarrow \Delta _X$).

Next, we define a certain quotient $\Delta_U\twoheadrightarrow \Delta_U^{p,i+1}$ of $\Delta_U$, which dominates $\Delta_{U,i+1/p}$. 
Let $i\ge 0$, and $G$ a finite quotient of $\Delta _{U,i+1/p}$ which inserts in the following commutative diagram.  
$$
\CD
1 @>>> I_{U,i+1/p}  @>>>  \Delta _{U,i+1/p} @>>> \Delta _X @>>> 1\\
@. @VVV     @VVV   @VVV\\
1 @>>> I  @>>> G @>>> {G }^{\et} @>>> 1\\
\endCD
$$
Thus, the quotient $G$ corresponds to a finite Galois cover $X_1'\to X$ with Galois group $G$, which factorizes as
$X_1'\to {X_1'}^{\et}\to X$, where ${X_1'}^{\et}\to X$ is the maximal \'etale sub-cover with Galois group $G^{\et}$, and $X_1'\to {X_1'}^{\et}$ is 
a tamely ramified Galois cover with Galois group $I$.  
Moreover, $I$ is an $(i+1)$-th step solvable group whose successive abelian quotients are annihilated by $p$.
We will assume for the remaining discussion, without loss of generality, that $G$ is {\it not} a quotient of $\Delta _{U,i/p}$.

Let $s\in \{1,\cdots,n\}$, and $I_{x_s}\subset G$ an inertia subgroup associated to $x_s$. 
Then $I_{x_s}$ is cyclic of order $p^t$, with ${t\le i+1}$, as follows from the structure of $I$. Write
$\Delta_{U,1/p^{i+1}}\defeq \Delta _U/I_U(1/p^{i+1})$
for the maximal cuspidally abelian annihilated by $p^{i+1}$ quotient of $\Delta_U$ (with respect to the 
surjection $\Delta_U\twoheadrightarrow \Delta_X$).
The following claim follows immediately from the well-known structure of $\Delta _U$ (cf. [Grothendieck], Expos\'e X, Corollaire 3.10).

\proclaim {Claim 1.2.1} There exists a finite quotient $G'$  of $\Delta_{U,1/p^{i+1}}$ which inserts in the following commutative diagram (where the vertical maps are surjective)
$$
\CD
1 @>>> I_{U,1/p^{i+1}}  @>>>  \Delta _{U,1/p^{i+1}} @>>> \Delta _X @>>> 1\\
@. @VVV     @VVV   @VVV\\
1 @>>> I'  @>>> G' @>>> {G '}^{\et} @>>> 1\\
\endCD
$$
such that the following holds. The quotient $\Delta_{U,1/p^{i+1}}\twoheadrightarrow G'$
corresponds to a Galois cover $X_2'\to X$ with Galois group $G'$ which factorizes as
$X_2'\to {X_2'}^{\et}\to X$, where ${X_2'}^{\et}\to X$ is the maximal \'etale sub-cover with Galois group ${G'}^{\et}$, and $X_2'\to {X_2'}^{\et}$ is a tamely ramified cover
with Galois group $I'$ (an abelian group annihilated by $p^{i+1}$). Further,
for $s\in \{1,\cdots,n\}$, and $I'_{x_s}\subseteq I'$ an inertia subgroup associated to $x_s$,
then $I'_{x_s}$ is {\it cyclic} of order $p^{i+1}$. 
\endproclaim


Let $K_1'\defeq K_{X_1'}$ (resp. $K_2'\defeq K_{X_2'}$) be the function field of $X_1'$ (resp. $X_2'$),
$L'\defeq K_1'.K_2'$ the compositum of $K_1'$ and $K_2'$, and $Y$ the normalisation of $X$ in $L'$. Thus, $Y\to X$ is a Galois cover with 
Galois group $H\subseteq G\times G'$ which is \'etale above $U$.
Note that $H$ maps onto $G$, $G'$, and the quotient $\Delta _U\twoheadrightarrow H$ doesn't factorize through $\Delta _U\twoheadrightarrow \Delta _{U,i+1/p}$
if $i\ge 1$.

Let $I''\defeq \Ker (I\twoheadrightarrow I^{\ab})$. Thus, $I''$ is an $i$-th step solvable group whose successive quotients are annihilated by $p$,
and is a characteristic subgroup of $I$. 
Write $\widetilde G\defeq G/I''$, and let $\widetilde H$ be the image of $H$ in  the quotient $\widetilde G\times G'$ of $G\times G'$.  We have a commutative diagram of exact sequences 
$$
\CD
1  @>>> I_H\defeq H\cap \left (I''\times \{1\} \right) @>>> H @>>> \widetilde H @>>>1\\
@. @VVV  @VVV  @VVV  \\
1 @>>> I''\times \{1\}  @>>> G\times G' @>>> \widetilde G\times G' @>>> 1\\
\endCD
$$
where the vertical maps are natural inclusions.

\proclaim {Lemma 1.2.2}
The group $\widetilde H$ is a quotient of $\Delta _{U,1/p^{i+1}}$. Moreover, the Galois cover $Y\to X$ factorizes as $Y\to Y'\to X$,
where $Y'\to X$ is a tamely ramified Galois cover with Galois group $\widetilde H$, and $Y\to Y'$ is an {\bf \'etale} Galois cover with Galois group $I_H\subseteq I''$: an {\bf $i$-th step solvable} group whose successive abelian quotients are {\bf annihilated by $p$}.
\endproclaim

\demo{Proof} The first assertion follows from the fact that the inertia subgroup of $\widetilde H$ is a subgroup of $I^{\ab}\times I'$. The proof of the second assertion is 
similar to the proof of Lemma 1.1.4 using Abhyankar's Lemma (cf. loc. cit.).
\qed
\enddemo

The profinite group $\Delta _{U,1/p^{i+1}}$; being a quotient of $\Delta_U$, is topologically finitely generated. 
Hence there exists a sequence of characteristic open subgroups
$$...\subseteq \Delta _{U,1/p^{i+1}}[j+1]\subseteq \Delta _{U,1/p^{i+1}}[j]\subseteq...\subseteq \Delta _{U,1/p^{i+1}}[1]\defeq \Delta _{U,1/p^{i+1}}$$
of $\Delta_{U,1/p^{i+1}}$ such that 
$\bigcap _{j\ge 1}\Delta _{U,1/p^{i+1}}[j]=\{1\}$. 
The open subgroup $\Delta _{U,1/p^{i+1}}[j]\subseteq \Delta_{U,1/p^{i+1}}$ corresponds to a
finite Galois cover  $(X_{i+1,j}')^U\to X$  
between smooth connected and proper $\ell$-curves, with Galois group $G_{i+1,j}^U\defeq \Delta _{U,1/p^{i+1}}/\Delta _{U,1/p^{i+1}}[j]$,
and which restricts to an \'etale cover $(V_{i+1,j})^U\to U$. 
The geometric point $\eta$ determines a geometric point $\eta_{i+1,j}'$ of $(X_{i+1,j}')^U$ and $(V_{i+1,j})^U$. 
Write $(\Delta_{i+1,j}')^U=\Delta_{(X_{i+1,j}')^U}\defeq \pi_1((X_{i+1,j}')^U,\eta_{i+1,j}')^{\Sigma}$
for the maximal pro-$\Sigma$ \'etale fundamental group of $(X_{i+1,j}')^U$ with base point $\eta_{i+1,j}'$, and 
$((\Delta _{i+1,j}')^U)_{i/p}$
for the maximal $i/p$-th step prosolvable quotient of $(\Delta_{i+1,j}')^U$. Consider the following push-out diagram
$$
\CD
1 @>>> \pi_1((V_{i+1,j})^U,\eta_{i+1,j}')^{\Sigma}@>>>  \Delta _U @>>> G_{i+1,j}^U @>>> 1\\
@. @VVV     @VVV   @VVV\\
1 @>>> ((\Delta _{i+1,j}')^U)_{i/p}  @>>> \widetilde G_{i+1,j}^U  @>>> G _{i+1,j}^U @>>> 1\\
\endCD
$$
($\ker \lgroup \pi_1((V_{i+1,j}')^U,\eta_{i+1,j}')^{\Sigma}\twoheadrightarrow ((\Delta _{i+1,j}')^U)_{i/p}\rgroup$ is a normal subgroup of $\Delta_U$).

\proclaim {Lemma 1.2.3}
With the above notations,
$G$ is a quotient of $\widetilde G_{i+1,j}^U$ for some $j\ge 1$.
\endproclaim

\demo{Proof} Follows from Lemma 1.2.2 and the various definitions.
\qed
\enddemo

Let $$\Delta _U^{p,i+1}\defeq \underset{j\ge 1} \to{\varprojlim}\ \widetilde G_{i+1,j}^U,$$
where $\widetilde G_{i+1,j}^U$ is as in Lemma 1.2.3 (the $\{\widetilde G_{i+1,j}^U\}_{j\ge 1}$ form a projective system). 
Thus, it follows from the various definitions that we have a natural exact sequence
$$1\to \underset{j\ge 1} \to{\varprojlim}((\Delta _{i+1,j}')^U)_{i/p} \to \Delta _U^{p,i+1}\to \Delta _{U,1/p^{i+1}}\to 1,\tag {1.7}$$
where the $\{((\Delta _{i+1,j}')^U)_{i/p}\}_{j\ge 1}$ are defined as above.

\proclaim {Proposition 1.2.4} The profinite group $\Delta _{U,i+1/p}$ is a {\bf quotient} of $\Delta _U^{p,i+1}$.
\endproclaim

\demo {Proof}
Follows from the above discussion (cf. Lemma 1.2.3).
\qed
\enddemo

Similarly, write
$\Cal I_{X,i+1/p}\defeq \Cal I_X/\Cal I_X(i+1/p).$
Thus, $\Cal I_{X,i+1/p}$ is the maximal $i+1/p$-th step prosolvable quotient of $\Cal I_X$, and $\Cal I_{X,1/p}$ is the maximal abelian annihilated by $p$ quotient of 
$\Cal I_X$. 
Write 
$$G_{X,i+1/p}\defeq G_X/\Cal I_{X}(i+1/p).$$
We shall refer to $G _{X,i+1/p}$ 
(resp. $G_{X,1/p}$) as the maximal {\bf cuspidally $i+1/p$-th step prosolvable} (resp. maximal {\bf cuspidally abelian annihilated by $p$})
quotient of $G_X$ (with respect to the surjection $G_X\twoheadrightarrow \Delta _X$). Also, write
$G_{X,1/p^{i+1}}\defeq G_X/\Cal I_{X}(1/p^{i+1})$,
and
$G_X^{p,i+1}\defeq \underset{U} \to{\varprojlim} (\Delta_U^{p,i+1})$
where the limit is taken over all open subschemes $U\subseteq X$. 
We have the following exact sequence
$$1\to \underset{U} \to{\varprojlim}\lgroup \underset{j\ge 1} \to{\varprojlim}((\Delta _{i+1,j}')^U)_{i/p}\rgroup \to G_X^{p,i+1}\to G_{X,1/p^{i+1}}\to 1.\tag {1.8}$$

\proclaim {Lemma 1.2.5} The profinite group $G_{X,i+1/p}$ is a {\bf quotient} of $G_X^{p,i+1}$.
\endproclaim

\demo{Proof} Follows from the various definitions, and Proposition 1.2.4.
\qed
\enddemo


\comment
Write
$$I_{U,i/\Sigma'}\defeq \prod _{p\in \Sigma'} \lgroup I_U/I_U(i/p)\rgroup.$$
We will refer to $I_{U,i/\Sigma'}$ as the maximal {\bf $i/\Sigma'$-th step prosolvable} quotient of $I_U$, and $I_{U,1/\Sigma'}$ as the maximal {\bf abelian annihilated by $\Sigma'$} quotient of $I_U$. 
Note that there exists a natural exact sequence
$$ 1\to   I_U[i+1/\Sigma']\to I_{U,i+1/\Sigma'}\to I_{U,i/\Sigma' }\to 1 \tag 1.11$$
where $I_U[i+1/\Sigma']$ is the subgroup $I_U(i/\Sigma')/I_U(i+1/\Sigma')$ of $I_{U,i+1/\Sigma'}$. In particular, $I_U[i+1/\Sigma' ]$ is {\it abelian} product of groups annihilated by $p$, for $p\in \Sigma'$.


Write 
$$\Delta _{U,i/\Sigma'}\defeq \Delta_U/I_{U}(i/\Sigma')$$
for $i\ge 0$. (Recall $I_U(i/\Sigma')$ is a characteristic subgroup of $I_U$ hence is normal in $\Delta_U$.) We shall refer to $\Delta _{U,i/\Sigma'}$ 
(resp. $\Delta _{U,1/\Sigma'}$) as the {\bf maximal cuspidally $i/\Sigma'$-th step prosolvable} (resp. {\bf maximal cuspidally abelian annihilated by $\Sigma'$})
quotient of $\Delta _U$ (with respect to the projection $\Delta _U\twoheadrightarrow \Delta _X$).
We have the following commutative diagram of exact sequences

$$
\CD
 @. 1 @.  1@. \\
@. @VVV    @VVV\\
@. I_U[i+1/\Sigma']  @= I_U[i+1/\Sigma']\\
@.  @VVV   @VVV \\
1 @>>> I_{U,i+1/\Sigma'}  @>>>  \Delta _{U,i+1/\Sigma'} @>>> \Delta _X @>>> 1\\
@.  @VVV     @VVV   @VVV\\
1 @>>> I_{U,i/\Sigma'}  @>>>  \Delta _{U,i/\Sigma'} @>>> \Delta _X @>>> 1\\
@.  @VVV   @VVV \\
@. 1@. 1@.\\
\endCD
\tag 1.12$$

Let $i\ge 0$. The profinite group $\Delta _{U,i/\Sigma'}$, being a quotient of $\Delta_U$, is topologically {\it finitely generated} (recall $\char (\ell)\notin \Sigma$). 
Hence there exists a sequence of {\it characteristic} open subgroups
$$...\subseteq \Delta _{U,i/\Sigma'}[j+1]\subseteq \Delta _{U,i/\Sigma'}[j]\subseteq...\subseteq \Delta _{U,i/\Sigma'}[1]\defeq \Delta _{U,i/\Sigma'}$$
of $\Delta_{U,i/\Sigma'}$, where $j\ge 1$ ranges over all positive integers, such that 
$\bigcap _{j\ge 1}\Delta _{U,i/\Sigma'}[j]=\{1\}$. The open subgroup $\Delta _{U,i/\Sigma' }[j]\subseteq \Delta_{U,i/\Sigma'}$ corresponds to a (possibly ramified)
finite (Galois) cover  
$$\widetilde X_{i,j}\defeq \widetilde X_{i,j}^{\Sigma'}\to X$$   
between smooth and proper $\ell$-curves, which is \'etale above $U$. The geometric point $\eta$ determines a geometric point $\eta_{i,j}$ of $X_{i,j}$. Write
$$\Delta_{i,j}\defeq\Delta_{i,j}^{U}=\Delta_{\widetilde X_{i,j}}\defeq \pi_1(\widetilde X_{i,j},\eta_{i,j})^{\Sigma}$$
for the maximal pro-$\Sigma$ \'etale fundamental group of $\widetilde X_{i,j}$ with base point $\eta_{i,j}$, and 
$\Delta _{i,j}^{\ab}\defeq (\Delta _{i,j}^{U})^{\ab}$ for the maximal abelian quotient of 
$\Delta_{i,j}$. The following provides a description of the structure of the profinite group $I_U[i+1/\Sigma']$, for $i\ge 1$.

\proclaim {Proposition 1.3.1} Assume $i\ge 1$. Then there exists a natural isomorphism $I_U[i+1/\Sigma']\isom \underset{j\ge 1} \to{\varprojlim}\ 
 \lgroup \prod _{p\in \Sigma'} \lgroup \Delta_{i,j}^{\ab}/p\Delta_{i,j}^{\ab}\rgroup \rgroup\defeq  \underset{j\ge 1} \to {\varprojlim}\ \lgroup \prod _{p\in \Sigma'} \lgroup (\Delta_{i,j}^U)^{\ab}/p(\Delta_{i,j}^U)^{\ab}\rgroup \rgroup$.
\endproclaim

\demo {Proof} Follows from the various definitions and Proposition 1.2.1.
\qed
\enddemo

Similarly, for $i\ge 0$ write
$$I_{X,i/\Sigma'}\defeq \prod _{p\in \Sigma'} \lgroup I_X/I_X(i/p)\rgroup .$$
We will refer to $I_{X,i/\Sigma'}$ as the maximal {\bf $i/\Sigma'$-th step prosolvable} quotient of $I_X$, and $I_{X,1/\Sigma'}$ as the maximal {\bf abelian annihilated by $\Sigma'$} quotient of $I_X$. 
Note that there exists a natural exact sequence
$$ 1\to   I_X[i+1/\Sigma']\to I_{X,i+1/\Sigma'}\to I_{X,i/p}\to 1 \tag 1.13$$
where $I_X[i+1/\Sigma']$ is the subgroup $I_X(i/\Sigma')/I_X(i+1/\Sigma')$ of $I_{X,i+1/\Sigma'}$. In particular, $I_X[i+1/\Sigma']$ is {\it abelian} 
product of groups annihilated by $p$, for all $p\in \Sigma'$. Write 
$$G_{X,i/\Sigma'}\defeq G_X/I_{X}(i/\Sigma')$$
for $i\ge 0$. (Recall $I_X(i/\Sigma')$ is a characteristic subgroup of $I_X$ hence is normal in $G_X$.) We shall refer to $G _{X,i/\Sigma'}$ 
(resp. $G_{X,1/\Sigma'}$) as the {\bf maximal cuspidally $i/\Sigma'$-th step prosolvable} (resp. {\bf maximal cuspidally abelian annihilated by $\Sigma'$})
quotient of $G_X$ (with respect to the projection $G_X\twoheadrightarrow \Delta _X$).
We have the following commutative diagram of exact sequences

$$
\CD
 @. 1 @.  1@. \\
@. @VVV    @VVV\\
@. I_X[i+1/\Sigma']  @= I_X[i+1/\Sigma']\\
@.  @VVV   @VVV \\
1 @>>> I_{X,i+1/\Sigma'}  @>>>  G_{X,i+1/\Sigma'} @>>> \Delta _X @>>> 1\\
@.  @VVV     @VVV   @VVV\\
1 @>>> I_{X,i/\Sigma'}  @>>>  G _{X,i/\Sigma'} @>>> \Delta _X @>>> 1\\
@.  @VVV   @VVV \\
@. 1@. 1@.\\
\endCD
\tag 1.13$$

\proclaim {Lemma 1.3.2} There are natural isomorphisms $G_{X,i/\Sigma'}\isom \underset{U} \to{\varprojlim}\ \Delta _{U,i/\Sigma'}$, $I_{X,i/\Sigma'}\isom \underset{U} \to{\varprojlim}\ I _{U,i/\Sigma'}$,
and $I_{X}[i/\Sigma']\isom \underset{U} \to{\varprojlim}\ I _{U}[i/\Sigma']$, where the projective limit is taken over all open subschemes $U\subseteq X$. Moreover, for $i\ge 1$
we have a natural isomorphism
$I_{X}[i+1/\Sigma']\isom \underset{U} \to{\varprojlim}     \lgroup \underset{j\ge 1} \to{\varprojlim} \lgroup \prod _{p\in \Sigma'}
\lgroup {\Delta _{i,j}^{U}})^{\ab}/p {\Delta _{i,j}^{U}})^{\ab}\rgroup \rgroup$ 
(cf. discussion before Proposition 1.3.1 for the definition of 
$(\Delta _{i,j}^{U})^{\ab})$.
\endproclaim

\demo{Proof} Follows from the various definitions and Proposition 1.2.2.
\qed
\enddemo
\endcomment

\subhead
\S2. Cuspidalisation of sections of cuspidally $i$-th step prosolvable arithmetric  fundamental groups
\endsubhead
In this section $k$ is a field with $\char (k)=l\ge 0$, $X$ is a proper smooth and geometrically connected
hyperbolic (i.e., $\genus(X)\ge 2$) curve over $k$, and $K_X$ its function field.
Let $\eta $ be a geometric point of $X$ above
its generic point; it determines an algebraic closure $\overline  k$
of $k$, 
and a geometric point $\overline \eta$ of $\overline X\defeq X\times _k \overline k$.

\subhead 2.1
\endsubhead
Let $\Sigma \subseteq \Primes$
be a non-empty set of prime integers.
In case $\char (k)=l >0$ we assume that $l\notin \Sigma$. 
Write
$\Delta_X\defeq \pi_1(\overline X,\overline \eta)^{\Sigma}$
for the maximal pro-$\Sigma$ quotient of $\pi_1(\overline X,\overline \eta)$, and
$\Pi_X\defeq  \pi_1(X, \eta)/ \Ker  (\pi_1(\overline X,\overline \eta)\twoheadrightarrow
\pi_1(\overline X,\overline \eta)^{\Sigma})$. Thus, we have an exact sequence 
$$1\to \Delta_X\to \Pi_X @>{\pr_{X,\Sigma}}>> G_k\defeq \Gal(\overline k/k)\to 1.\tag {$2.1$}$$
We shall refer to 
$\pi_1(X, \eta)^{(\Sigma)}\defeq \Pi_X$
as the {\bf geometrically pro-$\Sigma$ arithmetic fundamental group} of $X$. 

\subhead {2.1.1}
\endsubhead
Let $U\subseteq X$ be a nonempty open subscheme. 
Write $\Delta_U\defeq \pi_1(\overline U,\overline \eta)^{\Sigma}$
for the maximal pro-$\Sigma$ quotient of the fundamental group $\pi_1(\overline U,\overline \eta)$ of $\overline U$
with base point $\overline \eta$, and 
$\Pi_U\defeq  \pi_1(U, \eta)/ \Ker  (\pi_1(\overline U,\overline \eta)\twoheadrightarrow
\pi_1(\overline U,\overline \eta)^{\Sigma})$.
Thus, we have an exact sequence
$1 \to \Delta_U  \to \Pi_U \to G_k \to 1$.
Let
$I_U\defeq \Ker (\Pi_U\twoheadrightarrow \Pi_X)=\Ker (\Delta_U\twoheadrightarrow \Delta_X)$
be the cuspidal subgroup of $\Pi_U$ (with respect to the surjection $\Pi_U\twoheadrightarrow \Pi_X$). 
We have the following exact sequence
$$1\to I_U\to \Pi_U\to \Pi_X\to 1.\tag {2.2}$$
Let $i\ge 0$ be an integer, $I_{U,i}$ the maximal $i$-th step prosolvable quotient of $I_U$ (cf. 1.1), and $\Pi_U^{(i-\sol)}\defeq \Pi_U/\Ker (I_U\twoheadrightarrow I_{U,i})$.
We shall refer to $\Pi_U^{(i-\sol)}$ as the maximal (geometrically)  {\bf cuspidally $i$-th step prosolvable quotient} of $\Pi_U$
(with respect to the surjection $\Pi_U\twoheadrightarrow \Pi_X$).

\subhead {2.1.2}
\endsubhead
Similarly, we have an exact sequence of absolute Galois groups
$1\to G_{\overline k.K_X}\to G_{K_X}\to G_k\to 1$,
where $G_{\overline k.K_X}\defeq \pi_1(\Spec (\overline k.K_X),\eta)$, and $G_{K_X}\defeq \pi_1(\Spec (K_X),\eta)$.
Let 
$G_{\overline X}\defeq G_{\overline k.K_X}^{\Sigma}$
be the maximal pro-$\Sigma$ quotient of $G_{\overline k.K_X}$,
$G_X\defeq  G_{K_X}/ \Ker  ( G_{\overline k.K_X}  \twoheadrightarrow G_{\overline k.K_X}^{\Sigma})$, and
$\Cal I_X\defeq \Ker (G_{X}\twoheadrightarrow \Pi_X)=\Ker (G_{\overline X}\twoheadrightarrow \Delta_X)$
the cuspidal subgroup of $G_{X}$ (with respect to the surjection $G_X\twoheadrightarrow \Pi_X$). 
Let $\Cal I_{X,i}$ be the maximal $i$-th step prosolvable quotient of $\Cal I_X$. By pushing the exact sequence
$1\to \Cal I_X\to G_{X}  \to \Pi_X\to 1$
by the characteristic quotient $\Cal I_X\twoheadrightarrow  \Cal I_{X,i}$ we obtain an exact sequence
$$1\to \Cal I_{X,i}\to G_{X}^{(i-\sol)}\to \Pi_X\to 1.\tag {$2.3$}$$
We will refer to the quotient $G_{X} ^{(i-\sol)}$ as the maximal (geometrically) {\bf cuspidally $i$-th step prosolvable quotient} of $ G_{X}$
(with respect to the surjective homomorphism $G_{X}\twoheadrightarrow \Pi_X$).
There exist natural isomorphisms
$$ G_{X} ^{(i-\sol)}\isom \underset{U}\to {\varprojlim}\ \Pi_U^{(i-\sol)},\ \ \ \ \ \ \ \ \Cal I_{X,i} \isom \underset{U}\to {\varprojlim}\ I_{U,i}\ ,$$ 
where the limit is taken over all open subschemes $U\subseteq X$.

\subhead
2.2
\endsubhead
Let $i\ge 0$. We have a commutative diagram of exact sequences
$$
\CD
@. 1  @. 1\\
@. @VVV  @VVV\\
@. \Cal I_{X}[i+1] @= \Cal I_{X}[i+1]\\
@. @VVV     @VVV\\
1 @>>> G_{\overline X,i+1}@>>>  G_X^{(i+1-\sol)} @>>> G_k@>>> 1\\
@.  @VVV  @VVV @VVV \\
1 @>>> G _{\overline X,i}  @>>> G_X^{(i-\sol)} @>>> G_k @>>> 1\\
@. @VVV   @VVV \\
@. 1 @. 1 @.\\
\endCD
\tag {2.4}
$$
and similarly, for $U\subseteq X$ nonempty open, we have the following commutative diagram
$$
\CD
@. 1  @. 1\\
@. @VVV  @VVV\\
@. I_{U}[i+1] @= I_{U}[i+1]\\
@. @VVV     @VVV\\
1 @>>> \Delta _{U,i+1}@>>>  \Pi_U^{(i+1-\sol)} @>>> G_k@>>> 1\\
@.  @VVV  @VVV @VVV \\
1 @>>> \Delta _{U,i}  @>>> \Pi_U^{(i-\sol)} @>>> G_k @>>> 1\\
@. @VVV   @VVV \\
@. 1 @. 1 @. \\
\endCD
\tag {2.5}
$$
(recall the definition of $\Cal I_{X}[i+1]$ and  $I_{U}[i+1]$ from $\S1$).

Assume that the lower horizontal sequence in diagram (2.4) splits. Let
$s:G_k\to G_X^{(i-\sol)}$
be a section of the projection $G_X^{(i-\sol)}\twoheadrightarrow G_k$, which induces a section
$s_U:G_k\to \Pi_U^{(i-\sol)}$
of the projection $\Pi_U^{(i-\sol)}\twoheadrightarrow G_k$, $\forall U\subseteq X$ open.

\proclaim {The cuspidalisation problem for sections of cuspidally $i$-th step prosolvable arithmetric  fundamental groups}
Let $i\ge 0$. Given a section  $s_U:G_k\to \Pi_U^{(i-\sol)}$ as above, is it possible to construct a section $\tilde s_U:G_k\to \Pi_U^{(i+1-\sol)}$
of the projection $\Pi_U^{(i+1-\sol)}\twoheadrightarrow G_k$ which {\bf lifts} the section $s_U$, i.e., which fits in a commutative diagram
$$
\CD
G_k @>{\tilde s_U}>> \Pi_U^{(i+1-\sol)}\\
@V{\id}VV   @VVV\\ 
G_k @>{s_U}>> \Pi_U^{(i-\sol)}\\
\endCD
$$
where the right vertical map is the natural surjection?

Similarly, is it possible to construct a section $\tilde s:G_k\to G_X^{(i+1-\sol)}$
of the projection $G_X^{(i+1-\sol)}\twoheadrightarrow G_k$ which {\bf lifts} the section $s$, i.e., which fits in a commutative diagram
$$
\CD
G_k @>{\tilde s}>> G_X^{(i+1-\sol)}\\
@V{\id}VV   @VVV\\ 
G_k @>{s}>> G_X^{(i-\sol)}\\
\endCD
$$
where the right vertical map is the natural surjection?
\endproclaim

\subhead 2.3
\endsubhead The above cuspidalisation problem in the case $i=0$ has been investigated in [Sa\"\i di].
Next, we will investigate this problem in the case ${i\ge 1}$. 

We use the notations in 2.2.
Recall the definition of the characteristic open subgroups $\{\Delta_{U,i}[j]\}_{j\ge 1}$
such that $\bigcap _{j\ge 1}\Delta _{U,i}[j]=\{1\}$ (cf. discussion before Proposition 1.1.1).  
Write ${\widehat \Pi _U}[i,j] \defeq {\widehat \Pi _U}[i,j][s_U] \defeq \Delta _{U,i}[j]. s_U (G_k).$
Thus, ${\widehat \Pi_U}[i,j]\subseteq \Pi_U^{(i-\sol)}$ is an open subgroup which contains the image 
$s_U(G_k)$ of the section $s_U$. Write 
${\Pi_U}[i,j]\defeq {\Pi_U}[i,j][s_U]$
for the inverse image of $\widehat \Pi _U[i,j]$
in $\Pi_U$. Thus, $\Pi_U[i,j]\subseteq \Pi_U$ is an open subgroup which corresponds
to an \'etale cover $V_{i,j}\to U$,
where $V_{i,j}$ is a geometrically irreducible $k$-curve (since $\Pi_U[i,j]$ maps onto $G_k$ via the natural projection $\Pi_U\twoheadrightarrow G_k$ by the very definition of $\Pi_U[i,j]$).  

Write $X^U_{i,j}$ (resp. $\overline X^U_{i,j}$) for the smooth compactification of $V_{i,j}$ (resp. $\overline V_{i,j}\defeq V_{i,j}\times _k\overline k$).
We have an exact sequence 
$1 \to  \pi_1(\overline X^U_{i,j},\overline \eta_{i,j})\to \pi_1(X^U_{i,j},\eta_{i,j}) \to  G_k \to 1$,
where $\eta_{i,j}$ (resp. $\overline \eta_{i,j}$) is a geometric point naturally induced by $\eta$. 
Write $\pi_1(\overline X^U_{i,j},\overline \eta_{i,j}) ^{\Sigma,\ab}$ for the maximal pro-$\Sigma$ abelian quotient of $ \pi_1(\overline X^U_{i,j},\overline \eta_{i,j})$,  and 
$\pi_1(X^U_{i,j},\eta_{i,j})^{(\Sigma,\ab)}\defeq \pi_1(X^U_{i,j},\eta_{i,j}) /\Ker ( \pi_1(\overline X^U_{i,j},\overline \eta_{i,j})\twoheadrightarrow  
\pi_1(\overline X^U_{i,j},\overline \eta_{i,j})^{\Sigma,\ab})$ for the geometrically pro-$\Sigma$ abelian fundamental group of $X^U_{i,j}$. 
Consider the following pull-back diagram.
$$
\CD
1 @>>>  I_{U}[i+1]   @>>>  \Cal H_{U,i}\defeq \Cal H_{U,i}[s_U] @>>>  G_k  @>>> 1\\
@.    @VVV     @VVV     @V{s_U}VV \\
1 @>>> I_{U}[i+1]  @>>> \Pi_U^{(i+1-\sol)} @>>> \Pi_U^{(i-\sol)} @>>> 1\\
\endCD
\tag 2.6
$$

\proclaim {Lemma 2.3.1} The section $s_U:G_k\to \Pi_U^{(i-\sol)}$ {\bf lifts} to a section $\tilde s_U:G_k\to \Pi_U^{(i+1-\sol)}$
of the projection $\Pi_U^{(i+1-\sol)}\twoheadrightarrow G_k$ {\bf if and only if} the group extension
$1 \to I_{U}[i+1]  \to  \Cal H_{U,i} \to  G_k  \to 1$ {\bf splits}.
\endproclaim

\demo{Proof}
Follows immediately from the diagram (2.6).
\qed
\enddemo

\proclaim {Lemma 2.3.2} Assume ${i\ge 1}$. 
Then we have natural identifications $I_U[i+1]  \isom \underset{j\ge 1} \to{\varprojlim}\  \pi_1(\overline X^U_{i,j},\overline \eta_{i,j})^{\Sigma,\ab}$,
and
$\Cal H_{U,i}\isom \underset{j\ge 1} \to{\varprojlim}\  \pi_1(X^U_{i,j},\eta_{i,j})^{(\Sigma,\ab)}$.
\endproclaim

\demo{Proof} Follows from the various definitions and Proposition 1.1.1.
\qed
\enddemo

\definition {Definition 2.3.3}  We say that the field $k$ satisfies the condition $\bold {(H_{\Sigma})}$ if the following holds. The Galois cohomology groups $H^1(G_k,M)$
are {\it finite} for every finite $G_k$-module $M$ whose cardinality is divisible only by primes in $\Sigma$. 
\enddefinition

\proclaim {Lemma 2.3.4} Assume that ${i\ge 1}$, and $k$ satisfies the condition $\bold {(H_{\Sigma})}$. Then the group extension 
$1\to I_{U}[i+1]  \to \Cal H_{U,i}\to  G_k  \to 1$ {\bf splits
if and only if} the group extensions $1\to  \pi_1(\overline X^U_{i,j},\overline \eta_{i,j})^{\Sigma,\ab} \to \pi_1(X^U_{i,j},\eta_{i,j})^{(\Sigma,\ab)} \to G_k \to 1$
{\bf split}, $\forall j\ge 1$.
\endproclaim

\demo{Proof}
The only if part follows immediately from Lemma 2.3.2. Conversely, assume that the group extension 
$\pi_1(X^U_{i,j},\eta_{i,j})^{(\Sigma,\ab)}$ splits, $\forall j\ge 1$.
Write $\Cal H_{U,i}=\underset{G} \to{\varprojlim}\  G$ as the projective limit of finite quotients $G$ which insert into a commutative diagram
$$
\CD
1 @>>>  I_{U}[i+1]   @>>>  \Cal H_{U,i}\defeq \Cal H_{U,i}[s_U] @>>>  G_k  @>>> 1\\
@.    @VVV     @VVV     @VVV \\
1 @>>> \overline G  @>>> G @>>> H @>>> 1\\
\endCD
$$
where the vertical maps are surjective. Let $\widetilde G$ be the pull-back of the group extension $G$ by the surjective homomorphism $G_k\twoheadrightarrow H$.
Thus, we have an exact sequence $1\to \overline G\to \widetilde G\to G_k\to 1$, and $\Cal H_{U,i}=\underset{\widetilde G} \to{\varprojlim}\  \widetilde G$. 
Given a (geometrically finite) quotient $\Cal H_{U,i}\twoheadrightarrow \widetilde G$ as above, it factorizes as $\Cal H_{U,i}\twoheadrightarrow 
\pi_1(X^U_{i,j},\eta_{i,j})^{(\Sigma,\ab)} \twoheadrightarrow \widetilde G$ for some $j\ge 1$ (cf. Lemma 2.3.2). In particular, the group extension $\widetilde G$
splits by our assumption that the group extensions $\pi_1(X^U_{i,j},\eta_{i,j})^{(\Sigma,\ab)}$ split, $\forall j\ge 1$.
The set $\Sect (G_k,\Cal H_{U,i})$ of continuous splittings 
of the group extension $\Cal H_{U,i}$ is naturally identified with the inverse limit $\underset{\widetilde G}\to {\varprojlim} \Sect (G_k,\widetilde G)$ of the sets of continuous splittings of the group extensions $\widetilde G$ as above. For a (geometrically finite) quotient $\widetilde G$ of $\Cal H_{U,i}$ as above
the set $\Sect (G_k,\widetilde G)$ is non-empty (cf. above discussion) and is, up to conjugation by elements of  $\overline G$, a torsor under the Galois cohomology group $H^1(G_k,\overline G)$, which is finite by our assumption that $k$ satisfies the condition $\bold {(H_{\Sigma})}$. 
Thus, the set $\Sect (G_k,\widetilde G)$ is a non-empty finite set. Hence the set
$\Sect (G_k,\Cal H_{U,i})$ is non-empty being the projective limit of non-empty finite sets.
\qed
\enddemo

For $j\ge 1$, let $J_{i,j}^U\defeq \Pic^0_k(X_{i,j}^U)$ be the jacobian of $X_{i,j}^U$, and  $(J_{i,j}^1)^U\defeq \Pic^1_k(X_{i,j}^U)$. 
Thus, $(J_{i,j}^1)^U$ is a torsor under $J_{i,j}^U$.

\proclaim {Lemma 2.3.5} The group extension 
$1\to  \pi_1(\overline X^U_{i,j},\overline \eta_{i,j})^{\Sigma,\ab} \to \pi_1(X^U_{i,j},\eta_{i,j})^{(\Sigma,\ab)} \to G_k \to 1$
{\bf splits if and only if} the class of $(J_{i,j}^1)^U$ in $H^1(G_k,J_{i,j}^U)$ lies in the maximal $\Sigma$-divisible subgroup of  $H^1(G_k,J_{i,j}^U)$, 
i.e., the maximal subgroup of $H^1(G_k,J_{i,j}^U)$ which is divisible by integers whose prime factors are in $\Sigma$.
\endproclaim

\demo{Proof}
This is a well known fact (cf. [Harari-Szamuely], Theorem 1.2,  and Proposition 2.1).
\qed
\enddemo

\proclaim{Proposition 2.3.6} We use the same notations as above. Assume that 
${i\ge 1}$, and $k$ satisfies the condition $\bold {(H_{\Sigma})}$. 
Then the section $s_U:G_k\to \Pi_U^{(i-\sol)}$ {\bf lifts} to a section $\tilde s_U:G_k\to \Pi_U^{(i+1-\sol)}$
of the projection $\Pi_U^{(i+1-\sol)}\twoheadrightarrow G_k$ {\bf if and only if} the class of $(J_{i,j}^1)^U$ in $H^1(G_k,J_{i,j}^U)$ 
lies in the maximal $\Sigma$-divisible subgroup of  $H^1(G_k,J_{i,j}^U)$, $\forall j\ge 1$.
\endproclaim

\demo{Proof} Follows from Lemmas 2.3.4 and 2.3.5.
\qed
\enddemo

Recall the section $s:G_k\to G_X^{(i-\sol)}$
of the projection $G_X^{(i-\sol)}\twoheadrightarrow G_k$ which induces the section
$s_U:G_k\to \Pi_U^{(i-\sol)}$ of $\Pi_U^{(i-\sol)}\twoheadrightarrow G_k$, $\forall U\subseteq X$ open (cf. 2.2). 

\proclaim {Lemma 2.3.7} Assume that $k$ satisfies the condition $\bold {(H_{\Sigma})}$. 
Then the section $s:G_k\to G_X^{(i-\sol)}$ {\bf lifts} to a section $\tilde s:G_k\to G_X^{(i+1-\sol)}$
of the projection $G_X^{(i+1-\sol)}\twoheadrightarrow G_k$ {\bf if and only if} for each nonempty open 
subscheme $U\subseteq X$ the section $s_U:G_k\to \Pi_U^{(i-\sol)}$ lifts to a section $\tilde s_U:G_k\to \Pi_U^{(i+1-\sol)}$
of the projection $\Pi_U^{(i+1-\sol)}\twoheadrightarrow G_k$.
\endproclaim

\demo{Proof} Similar to the proof of Lemma 2.3.4, using the facts that $k$ satisfies the condition $\bold {(H_{\Sigma})}$,
and $G_{X} ^{(i-\sol)}\isom \underset{U}\to {\varprojlim}\ \Pi_U^{(i-\sol)}$.
\qed
\enddemo

The following is our main result in this section.
\proclaim{Theorem 2.3.8} We use the same notations as above. Assume that  ${i\ge 1}$,
and $k$ satisfies the condition $\bold {(H_{\Sigma})}$. 
Then the section $s:G_k\to G_X^{(i-\sol)}$ {\bf lifts} to a section $\tilde s : G_k\to G_X^{(i+1-\sol)}$
of the projection $G_X^{(i+1-\sol)}\twoheadrightarrow G_k$ {\bf if and only if} $\forall U\subseteq X$ nonempty open subscheme the class of $(J_{i,j}^1)^U$ in $H^1(G_k,J_{i,j}^U)$ 
lies in the maximal $\Sigma$-divisible subgroup of  $H^1(G_k,J_{i,j}^U)$, $\forall j\ge 1$.
\endproclaim

\demo{Proof} Follows from Lemmas 2.3.4, 2.3.5, and 2.3.7.
\qed
\enddemo

\subhead
\S 3. Lifting of sections to cuspidally $2/p$-th step prosolvable arithmetric fundamental groups
\endsubhead
In this section we investigate a certain mod-$p$ variant of the cuspidalisation problem investigated in $\S2$ (as well as in [Sa\"\i di]). 
Throughout $\S3$ we use the same notations as in $\S2$.
Let $p\in \Sigma$ be a prime integer.

\subhead {3.1}
\endsubhead
Let $U\subseteq X$ be a nonempty open subscheme.
Let $i\ge 0$, $t\ge 1$, be integers, and $I_{U,i/p^t}$ the maximal $i/p^t$-th step prosolvable quotient of $I_U$ (cf. 1.2). By pushing the exact sequence $(2.2)$ by the 
surjective homomorphism $I_U\twoheadrightarrow I_{U,i/p^t}$ we obtain an exact sequence
$1 \to I_{U,i/p^t} \to \Pi_U^{(i/p^t-\sol)} \to  \Pi_X \to 1 .$
We shall refer to $\Pi_U^{(i/p^t-\sol)}$ as the maximal (geometrically)  {\bf cuspidally $i/p^t$-th step prosolvable quotient} of $\Pi_U$
(with respect to the surjection $\Pi_U\twoheadrightarrow \Pi_X$).
We have a commutative diagram of exact sequence.
$$
\CD
@. 1  @. 1\\
@. @VVV  @VVV\\
@. I_{U,i/p^t} @= I_{U,i/p^t}\\
@. @VVV     @VVV\\
1 @>>> \Delta _{U,i/p^t}@>>>  \Pi_U^{(i/p^t-\sol)} @>>> G_k@>>> 1\\
@.  @VVV  @VVV @VVV \\
1 @>>> \Delta _X  @>>> \Pi_X @>>> G_k @>>> 1\\
@. @VVV   @VVV \\
@. 1 @. 1 \\
\endCD
\tag {3.1}$$

Similarly, by pushing the exact sequence
$1\to \Cal I_X\to G_{X}  \to \Pi_X\to 1$ by the 
surjective homomorphism $\Cal I_X\twoheadrightarrow  \Cal I_{X,i/p^t}$ we obtain an exact sequence
$1\to \Cal I_{X,i/p^t}\to G_{X}^{(i/p^t-\sol)}\to \Pi_X\to 1.$
We will refer to $G_{X} ^{(i/p^t-\sol)}$ as the maximal (geometrically) {\bf cuspidally $i/p^t$-th step prosolvable quotient} of $ G_{X}$
(with respect to the surjective homomorphism $G_{X}\twoheadrightarrow \Pi_X$).
There exist natural isomorphisms
$$ G_{X} ^{(i/p^t-\sol)}\isom \underset{U}\to {\varprojlim}\ \Pi_U^{(i/p^t-\sol)},\ \ \ \ \ \ \ \ \Cal I_{X,i/p^t} \isom \underset{U}\to {\varprojlim}\ I_{U,i/p^t},$$ 
where the limits are over all open subschemes $U\subseteq X$, and a commutative diagram.
$$
\CD
@. 1  @. 1\\
@. @VVV  @VVV\\
@. \Cal I_{X,i/p^t} @= \Cal I_{X,i/p^t}\\
@. @VVV     @VVV\\
1 @>>> G_{\overline X,i/p^t}@>>>  G_X^{(i/p^t-\sol)} @>>> G_k@>>> 1\\
@.  @VVV  @VVV @VVV \\
1 @>>> \Delta _X  @>>> \Pi_X @>>> G_k @>>> 1\\
@. @VVV   @VVV \\
@. 1 @. 1 \\
\endCD
\tag {3.2}$$

\subhead
3.2
\endsubhead
Assume that the lower horizontal exact sequence in diagram (3.2) splits. Let
$s:G_k\to \Pi_X$
be a section of the projection $\Pi_X\twoheadrightarrow G_k$.

\proclaim {The lifting problem to sections of cuspidally $i+1/p^t$-th step prosolvable arithmetric  fundamental groups}
Let $i\ge 0$, $t\ge 1$, be integers. Given a section  $s:G_k\to \Pi_X$ as above is it possible to construct a section $s_{U,i+1}:G_k\to \Pi_U^{(i+1/p^t-\sol)}$
of the projection $\Pi_U^{(i+1/p^t-\sol)}\twoheadrightarrow G_k$ which {\bf lifts} the section $s$, i.e., which fits in a commutative diagram
$$
\CD
G_k @>{s_{U,i+1}}>> \Pi_U^{(i+1/p^t-\sol)}\\
@V{\id}VV   @VVV\\ 
G_k @>{s}>> \Pi_X\\
\endCD
$$
where the right vertical map is the natural surjection? 

Similarly, is it possible to construct a section $s_{i+1}:G_k\to G_X^{(i+1/p^t-\sol)}$
of the projection $G_X^{(i+1/p^t-\sol)}\twoheadrightarrow G_k$ which {\bf lifts} the section $s$, i.e., which fits in a commutative diagram
$$
\CD
G_k @>{s_{i+1}}>> G_X^{(i+1/p^t-\sol)}\\
@V{\id}VV   @VVV\\ 
G_k @>{s}>> \Pi_X\\
\endCD
$$
where the right vertical map is the natural surjection?
\endproclaim

\subhead 3.3.  The quotients $G_X\twoheadrightarrow G_X^{(p,i+1)}$, $\Pi_U\twoheadrightarrow \Pi_U^{(p,i+1)}$, and lifting of sections
\endsubhead
Next, recall the notations in 1.2 and the discussion therein, especially the definition of the quotient $\Delta _U\twoheadrightarrow \Delta_U^{p,i+1}$ (cf. the discussion after Lemma 1.2.3).  
The kernel of the surjective 
homomorphism $\Delta_U\twoheadrightarrow \Delta_U^{p,i+1}$
is a normal subgroup of $\Pi_U$ (as one easily verifies). 
Write $\Pi _U^{(p,i+1)}\defeq \Pi_U/\Ker (\Delta_U\twoheadrightarrow \Delta _U^{p,i+1}).$ Thus, we have an exact sequence
$$1\to  \Delta _U^{p,i+1}  \to  \Pi _U^{(p,i+1)} \to  G_k \to 1.\tag 3.3$$
Recall the exact sequence
$1\to \underset{j\ge 1} \to{\varprojlim}((\Delta _{i+1,j}')^U)_{i/p} \to \Delta _U^{p,i+1}\to \Delta _{U,1/p^{i+1}}\to 1$ (cf. loc. cit.).
The quotient $\Pi_U\twoheadrightarrow \Pi _{U}^{(1/p^{i+1}-\sol)}$ (cf. 3.1) factorizes through  $\Pi_U\twoheadrightarrow \Pi_U^{(p,i+1)}$ 
(cf. exact sequence (1.7)),
and we have 
a commutative diagram of exact sequences.
$$
\CD
@. 1 @. 1 @.\\
@. @VVV  @VVV\\
@. \underset{j\ge 1} \to{\varprojlim}((\Delta _{i+1,j}')^U)_{i/p} @= \underset{j\ge 1} \to{\varprojlim}((\Delta _{i+1,j}')^U)_{i/p} \\
@. @VVV  @VVV \\
1 @>>> \Delta _U^{p,i+1} @>>> \Pi_U^{(p,i+1)} @>>> G_k @>>> 1\\
@. @VVV @VVV @VVV\\
1@>>> \Delta _{U, 1/p^{i+1}} @>>>  \Pi _{U}^{(1/p^{i+1}-\sol)} @>>> G_k @>>> 1\\
@. @VVV  @VVV \\
@.  1  @. 1 \\
\endCD
\tag 3.4
$$

Similarly,  $\Ker (G_{\overline X}\twoheadrightarrow G_{\overline X}^{p,i+1})$ is a normal subgroup of $G_X$, and we have an exact sequence
$$1\to  G_{\overline X}^{p,i+1}  \to  G _X^{(p,i+1)} \to  G_k \to 1,\tag {3.5}$$
where 
$G _X^{(p,i+1)}\defeq G_X/  \Ker (G_{\overline X}\twoheadrightarrow G_{\overline X}^{p,i+1})$.
The exact sequence (1.8)
induces a commutative diagram of exact sequences.
$$
\CD
@. 1 @. 1 @.\\
@. @VVV  @VVV\\
@. \underset{U} \to{\varprojlim}\lgroup \underset{j\ge 1} \to{\varprojlim}((\Delta _{i+1,j}')^U)_{i/p}\rgroup @=
\underset{U} \to{\varprojlim}\lgroup \underset{j\ge 1} \to{\varprojlim}((\Delta _{i+1,j}')^U)_{i/p}\rgroup       \\
@. @VVV  @VVV \\
1 @>>> G_{\overline X}^{p,i+1} @>>> G_X^{(p,i+1)} @>>> G_k @>>> 1\\
@. @VVV @VVV @VVV\\
1@>>> G _{\overline X, 1/p^{i+1}} @>>>  G_{X}^{(1/p^{i+1}-\sol)} @>>> G_k @>>> 1\\
@. @VVV  @VVV \\
@.  1  @. 1 \\
\endCD
\tag 3.6
$$

Furthermore, $\Pi_U^{(i+1/p-\sol)}$ (resp. $G_X^{(i+1/p-\sol)}$) is a quotient of $\Pi_U^{(p,i+1)}$ (resp. of $G_X^{(p,i+1)}$) (cf. Lemmas 1.2.4 and 1.2.5)
and we have commutative diagrams
$$
\CD
1 @>>> \Delta _U^{p,i+1} @>>> \Pi_U^{(p,i+1)} @>>> G_k @>>> 1\\
@. @VVV @VVV @VVV\\
1 @>>> \Delta _{U,i+1/p} @>>>  \Pi_U^{(i+1/p-\sol)} @>>> G_k @>>> 1\\
\endCD
\tag 3.7
$$
 resp.
$$
\CD
1 @>>> G_{\overline X}^{p,i+1} @>>> G_X^{(p,i+1)} @>>> G_k @>>> 1\\
@. @VVV @VVV @VVV\\
1 @>>> G _{\overline X,i+1/p} @>>>  G_X^{(i+1/p-\sol)} @>>> G_k @>>> 1\\
\endCD
\tag 3.8
$$
where the left and middle vertical maps are surjective.
\proclaim {The lifting problem to sections of $\Pi_U^{(p,i+1)}$, and $G_X^{(p,i+1)}$}
Given a section  $s:G_k\to \Pi_X$ of the projection $\Pi_X\twoheadrightarrow G_k$ as in $3.2$, is it possible to construct a section $\tilde s_{U,i+1}:G_k\to \Pi_U^{(p,i+1)}$
of the projection $\Pi_U^{(p,i+1)}\twoheadrightarrow G_k$ which {\bf lifts} the section $s$, i.e., which fits in a commutative diagram
$$
\CD
G_k @>{\tilde s_{U,i+1}}>> \Pi_U^{(p,i+1)}\\
@V{\id}VV   @VVV\\ 
G_k @>{s}>> \Pi_X\\
\endCD
$$
where the right vertical map is the natural surjection?

Similarly, is it possible to construct a section $\tilde s_{i+1}:G_k\to G_X^{(p,i+1)}$
of the projection $G_X^{(p,i+1)}\twoheadrightarrow G_k$ which {\bf lifts} the section $s$, i.e., which fits in a commutative diagram
$$
\CD
G_k @>{\tilde s}_{i+1}>> G_X^{(p,i+1)}\\
@V{\id}VV   @VVV\\ 
G_k @>{s}>> \Pi_X\\
\endCD
$$
where the right vertical map is the natural surjection?
\endproclaim

\proclaim {Lemma 3.3.1} A positive answer to the lifting problem posed in 3.3 implies a positive answer to the lifting problem posed in 3.2
in the case ${t=1}$.
\endproclaim

\demo {Proof}
Follows immediately from the above commutative diagrams (3.7) and (3.8).
\qed
\enddemo

\subhead {3.4}
\endsubhead
In this section we investigate the lifting problem posed in 3.3 in the case $i=1$, and draw consequences for the lifting problem posed in 3.2 in the case ${t=i=1}$
(the only case we need for applications in $\S4$).
Let $s:G_k\to \Pi_X$
be a section of the projection $\Pi_X\twoheadrightarrow G_k$.

\subhead{3.4.I}
\endsubhead
First, we investigate the problem of lifting the section $s$ to a section 
$s_{U,i+1}':G_k\to \Pi_U^{(1/p^{i+1}-\sol)}$ (resp. 
$s_{i+1}':G_k\to G_X^{(1/p^{i+1}-\sol)}$)
of the projection  $\Pi_U^{(1/p^{i+1}-\sol)}\twoheadrightarrow G_k$ (resp. $G_X^{(1/p^{i+1}-\sol)}\twoheadrightarrow G_k$).
Recall the notations in 2.3, and ($\forall U\subseteq X$ open) the sequence of characteristic open subgroups
$$...\subseteq \Delta_X[j+1]\defeq \Delta _{U,0}[j+1]\subseteq \Delta_X[j]\defeq \Delta _{U,0}[j]\subseteq...\subseteq \Delta_{X}[1]\defeq \Delta _{U,0}[1]=\Delta _X$$
of $\Delta_{U,0}=\Delta_X$ with 
$\bigcap _{j\ge 1}\Delta _{X}[j]=\{1\}.$ 
The sequence of (geometrically characteristic) open subgroups
$$... \subseteq  \Pi_X[j+1]\defeq \Pi_{U,0}[j+1]\subseteq  \Pi_X[j]\defeq \Pi_{U,0}[j]\subseteq...\subseteq \Pi_X[1]\defeq \Pi_{U,0}[1]=\Pi _X,$$ 
where 
$\Pi_{X}[j]\defeq \Delta _X[j].s(G_k)$,
corresponds to a tower of finite (not necessarily Galois) \'etale covers
$$... \to X_{j+1}\defeq X_{0,j+1}^U\to X_{j}\defeq X_{0,j}^U\to ...\to X\defeq X_{0,1}^U.$$
Note that $\Pi_{X}[j]$ identifies naturally with 
$\Pi_{X_{j}}\defeq \pi_1(X_{j},\eta_j)^{(\Sigma)}$, where $\eta_j$ is the base point induced by $\eta$.
Moreover, the section $s$ restricts to a section 
$s:G_k\to \Pi_{X_j}$ 
of the projection $\Pi_{X_j}\twoheadrightarrow G_k$, $\forall j\ge 1$.

Let $i\ge 0$, $j\ge 1$,  be integers. Recall the Kummer sequence 
$$1\to \mu_{p^{i+1}}\to \Bbb G_m @>p^{i+1}>> \Bbb G_m\to 1$$
in \'etale topology, which induces an exact sequence 
$$0 \to \Pic (X_j)/p^{i+1} \Pic (X_j) \to  H^2(X_j,\mu_{p^{i+1}}) \to _{p^{i+1}}\Br (X_j)\to 0.$$
Here $\Pic\defeq H_{\et}^1(\ ,\Bbb G_m)$ is the Picard group, $\Br\defeq H_{\et}^2(\ ,\Bbb G_m)$
the Brauer-Grothendieck cohomological group, and $_{p^{i+1}} \Br\subseteq \Br$
the subgroup of $\Br$ which is annihilated by $p^{i+1}$. 
We identify $\Pic (X_j)/p^{i+1} \Pic (X_j)$ with its image in $H^2(X_j,\mu_{p^{i+1}})$
and refer to it as the Picard part of $H^2(X_j,\mu_{p^{i+1}})$.
By pulling back cohomology classes via the section $s:G_k\to \Pi_{X_j}$,
and bearing in mind the natural identification 
$H^2(\Pi_{X_j},\mu_{p^{i+1}})\isom H^2(X_j,\mu_{p^{i+1}})$
(cf. [Mochizuki], Proposition 1.1), we obtain a restriction homomorphism
$s_j^{\star}: H^2(X_j,\mu_{p^{i+1}})\to H^2(G_k,\mu_{p^{i+1}})$.

Observe that if $k'/k$ is a finite extension, and $X_{k'}\defeq X\times _kk'$,
then we have a cartesian diagram:
$$
\CD
1 @>>>  \Delta _{X_{k'}}     @>>> \Pi_{X_{k'}}   @>>> G_{k'} @>>> 1\\
  @.        @VVV           @VVV              @VVV \\
1 @>>>    \Delta _X      @>>> \Pi_X       @>>> G_k   @>>> 1
\endCD
$$
and the section $s$
induces a section
$s_{k'}:G_{k'}\to  \Pi_{X_{k'}}$ 
of the projection $\Pi_{X_{k'}}\twoheadrightarrow G_{k'}$.

\definition{Definition 3.4.1 (Sections with Cycle Classes Orthogonal to $\Pic$ mod-$p^{i+1}$)} (Compare with [Sa\"\i di], 1.4.)

\noindent
(i) We say that the section $s$ has a {\bf cycle class orthogonal to} $\Pic$ {\bf mod-$p^{i+1}$} if 
the homomorphism
$s_j^{\star}: H^2(X_j,\mu_{p^{i+1}})\to H^2(G_k,\mu_{p^{i+1}})$ 
annihilates the Picard part $\Pic (X_j)/p^{i+1}\Pic (X_j)$ of
$H^2(X_j,\mu_{p^{i+1}})$, $\forall j\ge 1$. 

\noindent
(ii)\ We say that the section $s$ has a {\bf cycle class uniformly orthogonal to $\Pic $ mod-$p^{i+1}$} 
(relative to the system of neighbourhoods $\{X_j\}_{j\ge 1}$ of $s$) if, for every finite extension $k'/k$, 
the induced section $s_{k'}:G_{k'}\to \Pi_{X_{k'}}$ has a cycle class 
orthogonal to $\Pic$ mod-$p^{i+1}$ (relative to the system of neighbourhoods of $s_{k'}$ which is induced by the $\{X_j\}_{j\ge 1}$).
\enddefinition

\definition {Definition 3.4.2}  We say that the field $k$ satisfies the condition $\bold {(H_{p^{i+1}})}$ if the following holds. The Galois cohomology groups $H^1(G_k,M)$
are {\it finite} for every finite $G_k$-module $M$ annihilated by $p^{i+1}$. 
\enddefinition

\proclaim {Theorem 3.4.3 (Lifting of Sections to Cuspidally mod-$p^{i+1}$ abelian Arithmetic Fundamental Groups)} 
Assume that $k$ satisfies the condition $\bold {(H_{p^{i+1}})}$ (cf. Definition 3.4.2). 
Let $s:G_k\to \Pi_X$ be a section of the projection $\Pi_X\twoheadrightarrow G_k$. Assume that $s$ has a {\bf cycle class uniformly orthogonal to $\Pic$ mod-$p^{i+1}$} (cf. Definition 3.4.1(ii)). 
Let $U\subseteq X$ be a nonempty open subscheme. Then there exists a section $s_{U,i+1}':G_k\to \Pi_{U}^{(1/p^{i+1}-\sol)}$
of the projection  $\Pi_{U}^{(1/p^{i+1}-\sol)}\twoheadrightarrow G_k$ which {\bf lifts} the section $s$, i.e., which inserts into
the following commutative diagram.
$$
\CD
G_k @>s_{U,i+1}'>>  \Pi_{U}^{(1/p^{i+1}-\sol)} \\
@V{\id}VV     @VVV  \\
G_k   @>{s}>> \Pi_X
\endCD
$$
\endproclaim

\demo{Proof} Similar to the proof of Theorem 2.3.3 in [Sa\"\i di].
\qed
\enddemo

\proclaim {Theorem 3.4.4 (Lifting of Sections to Cuspidally mod-$p^{i+1}$ abelian Galois 
Groups)} Assume that the field $k$ satisfies the condition $\bold {(H_{p^{i+1}})}$ (cf. Definition 3.4.2). 
Let $s:G_k\to \Pi_X$ be a section of the projection $\Pi_X\twoheadrightarrow G_k$. Then $s$ has a {\bf cycle class 
uniformly orthogonal to $\Pic$ mod-$p^{i+1}$} (cf. Definition 3.4.1 (ii))
{\bf if and only if} there exists a section $s_{i+1}':G_k\to G_{X}^{(1/p^{i+1}-\sol)}$ of the
projection  $G_{X}^{(1/p^{i+1}-\sol)}  \twoheadrightarrow G_k$ which {\bf lifts} the section $s$,
i.e., which inserts in the following commutative diagram.
$$
\CD
G_k @>s_{i+1}'>>   G_{X}^{(1/p^{i+1}-\sol)}\\
@V{\id}VV     @VVV  \\
G_k   @>{s}>> \Pi_X
\endCD
$$

\endproclaim

\demo{Proof} Similar to the proof of Theorem 2.3.5 in [Sa\"\i di].
\qed
\enddemo

\subhead{3.4.II}
\endsubhead
Next, let  
$s':G_k\to G_{X}^{(1/p^{2}-\sol)}$ 
be a section of the projection  $G_{X}^{(1/p^{2}-\sol)}  \twoheadrightarrow G_k$, which induces for every open subscheme $U\subseteq X$
a section 
$s_{U}':G_k\to \Pi_{U}^{(1/p^{2}-\sol)}$
of the projection  $\Pi_{U}^{(1/p^{2}-\sol)}\twoheadrightarrow G_k$. We investigate the problem of lifting the section
$s'_U$ (resp. $s'$) to a section $\tilde s_{U}:G_k\to \Pi_U^{(p,2)}$ (resp. $\tilde s:G_k\to G_X^{(p,2)}$)
of the projection $\Pi_U^{(p,2)}\twoheadrightarrow G_k$ (resp. $G_X^{(p,2)}\twoheadrightarrow G_k$)
(cf. diagrams (3.4) and (3.6)). Let $U\subseteq X$ be an open subscheme. Consider the following pull-back diagram.
$$
\CD
1 @>>>  \underset{j\ge 1} \to{\varprojlim}((\Delta _{2,j}')^U)_{1/p}@>>>  \Cal H_U^{(p,2)} 
@>>> G_k @>>>1\\
@. @VVV @VVV @V{s_{U}'}VV \\
1 @>>>  \underset{j\ge 1} \to{\varprojlim} ((\Delta _{2,j}')^U)_{1/p} @>>> \Pi_U^{(p,2)}@>>>  \Pi _{U}^{(1/p^{2}-\sol)} @>>> 1\\
\endCD
\tag 3.9
$$

\proclaim {Lemma 3.4.5} The section $s_{U}'$ lifts to a section $\tilde s_{U}:G_k\to \Pi_U^{(p,2)}$
of the projection $\Pi_U^{(p,2)}\twoheadrightarrow G_k$ if and only if the group extension $\Cal H_U^{(p,2)}$ splits. 
\endproclaim

\demo{Proof} Follows immediately from diagram (3.9).
\qed
\enddemo

Recall the discussion and notations after Lemma 1.2.2, especially the definition of the 
$\{\Delta _{U,1/p^{i+1}}[j] \}_{j\ge 1}$. For $j\ge 1$, write
$\Pi_{U,1/p^{2}}[j]\defeq \Delta _{U,1/p^{2}}[j].s_{U}'(G_k).$
Thus, $\Pi_{U,1/p^{2}}[j]\subseteq  \Pi _{U}^{(1/p^{2}-\sol)}$ is an open subgroup corresponding to a (possibly tamely ramified)
cover
$\widetilde X_{j}^U\to X$   
between smooth, proper, and geometrically connected $k$-curves.
The geometric point $\eta$ determines a geometric point $\eta_{j}$ of $\widetilde X_{j}^U$. Write
$\Pi_{j}^U=\Pi_{j}^U[s_{U}']
\defeq \Pi_{\widetilde X_{j}^U}\defeq \pi_1(\widetilde X_{j}^U,\eta_{j})^{(\Sigma)}$, which inserts in the exact sequence
$1\to  \Delta _{j}^U\to \Pi_{j}^U  \to G_k\to 1$,
where
$\Delta _{j}^U\defeq \Delta _{{\widetilde X_{j}^U}\times_k \overline k}.$
Further, consider the push-out diagram
$$
\CD
1@>>> \Delta _{j}^U@>>> \Pi_{j}^U @>>> G_k @>>> 1\\
@. @VVV @VVV @VVV\\
1@>>> (\Delta _{j}^U)_{1/p}@>>> (\Pi_{j}^U)^{(1/p-\sol)}@>>>  G_k @>>> 1\\
\endCD
$$
which defines the geometrically $1/p$-th step solvable quotient $\Pi_{j}^U\twoheadrightarrow (\Pi_{j}^U)^{(1/p-\sol)}$ of $\Pi_{j}^U$.

\proclaim {Lemma 3.4.6} There are natural isomorphisms
$\underset{j\ge 1} \to{\varprojlim} ((\Delta _{2,j}')^U)_{1/p} \isom \underset{j\ge 1} \to{\varprojlim} (\Delta _{j}^U)_{1/p}$,
and
$\Cal H_U^{(p,2)}\isom  \underset{j\ge 1} \to{\varprojlim}  (\Pi _{j}^U)^{(1/p-\sol)}$.
\endproclaim

\demo{Proof} Follows from the various definitions. 
\qed
\enddemo
Note that $(\Delta _{2,j}')^U=\Delta _{j}^U$, we will in the sequel
write $\Delta _{j}^U$ instead of $(\Delta _{2,j}')^U$.

\proclaim {Lemma 3.4.7} Assume that $k$ satisfies the condition $\bold {(H_{p})}$ (cf. Definition 3.4.2). 
Then the group extension $\Cal H_U^{(p,2)}$ splits {\bf if and only if} the group extension $(\Pi_{j}^U)^{(1/p-\sol)}$
splits, $\forall j\ge 1$.
\endproclaim

\demo{Proof} Similar to the proof of Lemma 2.3.4, using the fact that $H^1(G_k,(\Delta _{j}^U)_{1/p})$ is finite if $k$ satisfies $\bold {(H_p)}$.
\qed
\enddemo

For $j\ge 1$, let $J_{j}[U]\defeq \Pic^0_k(\widetilde X_{j}^U)$ be the jacobian of $\widetilde X_{j}^U$, and  
$J_{j}^1[U]\defeq \Pic^1_k(\widetilde X_{j}^U)$. 

\proclaim {Lemma 3.4.8} The group extension  $(\Pi_{j}^U)^{(1/p-\sol)}$ splits {\bf if and only if}
the class of $ J_{j}^1[U]$ in $H^1(G_k,J_{j}[U])$ is divisible by $p$.

\endproclaim

\demo{Proof} This fact is well-known, see [Harari-Szamuely] for instance. 
Strictly speaking loc. cit. treats the splittings of the group extension
$(\Pi_{j}^U)^{(\ab)}$= the geometrically abelian quotient of $\Pi_{j}^U$, but a similar argument leads to  a mod-$p$
variant as above for any prime $p\in \Sigma$. 
\qed
\enddemo

\proclaim {Theorem 3.4.9} With the above notations, assume that $k$ satisfies the condition $\bold {(H_{p})}$ (cf. Definition 3.4.2).  Then
the section $s_{U}'$ lifts to a section $\tilde s_{U}:G_k\to \Pi_U^{(p,2)}$
of the projection $\Pi_U^{(p,2)}\twoheadrightarrow G_k$ {\bf if and only if} the class of $ J_{j}^1[U]$ in 
$H^1(G_k,J_{j}[U])$ is divisible by $p$, $\forall j\ge 1$.
\endproclaim

\demo{Proof} Follows from Lemmas 3.4.5, 3.4.7, and 3.4.8.
\qed
\enddemo

\proclaim {Theorem 3.4.10} With the above notations, assume that $k$ satisfies the condition $\bold {(H_{p})}$ (cf. Definition 3.4.2). Then
the section $s'$ lifts to a section $\tilde s:G_k\to G_X^{(p,2)}$
of the projection $G_X^{(p,2)}\twoheadrightarrow G_k$ {\bf if and only if} the class of $ J_{j}^1[U]$ in 
$H^1(G_k,J_{j}[U])$ is divisible by $p$, $\forall j\ge 1$, and 
$\forall U\subseteq X$ nonempty open subscheme as in the above discussion. 
\endproclaim

\demo{Proof} Similar to the proof of Theorem 3.4.9.
\qed
\enddemo

The following is our main result in this section.

\proclaim{Theorem 3.4.11} With the above notations, assume that the field $k$ satisfies the condition $\bold {(H_{p^{2}})}$ (cf. Definition 3.4.2). 
Let $s:G_k\to \Pi_X$ be a section of the projection $\Pi_X\twoheadrightarrow G_k$. Then $s$ lifts to a section $\tilde s:G_k\to G_X^{(p,2)}$
(resp. $s':G_k\to G_X^{(2/p)}$)
of the projection $G_X^{(p,2)}\twoheadrightarrow G_k$ (resp. $G_X^{(2/p)}\twoheadrightarrow G_k$)
{\bf if and only if} (resp. {\bf if}) the following two conditions occur.

{(i)} The section $s$ has a {\bf cycle class 
uniformly orthogonal to $\Pic$ mod-$p^{2}$} (cf. Definition 3.4.1 (ii))

{(ii)} There exists a section $s':G_k\to G_{X}^{(1/p^{2}-\sol)}$ of the
projection  $G_{X}^{(1/p^{2}-\sol)}  \twoheadrightarrow G_k$ which {\bf lifts} the section $s$ (this holds if (i) holds by Theorem 3.4.4) such that 
the class of $ J_{j}^1[U]$ in 
$H^1(G_k,J_{j}[U])$ is {\bf divisible} by $p$, $\forall j\ge 1$, and 
$\forall U\subseteq X$ nonempty open subscheme.
\endproclaim

\demo{Proof} Follows from Theorems 3.4.4 and 3.4.10. The resp. assertion follows from diagram (3.8) (cf. Lemma 3.3.1).
\qed
\enddemo

\subhead
\S 4. Geometric sections of arithmetic fundamental groups of $p$-adic curves
\endsubhead

\noindent
In this section, applying the results in $\S3$, we provide a characterisation of sections of (geometrically pro-$\Sigma$, $p\in \Sigma$) 
arithmetic fundamental groups of $p$-adic curves which arise from rational points. 
We use the notations in $\S2$ and $\S3$. 

Let $p$ be a prime integer. In this section $k$ is a $p$-adic local field, i.e., $k/\Bbb Q_p$ is a finite extension,
and we assume $p\in \Sigma.$ 
Let $s:G_k\to \Pi_X$
be a section of the projection $\Pi_X\twoheadrightarrow G_k$.

\definition{Definition 4.1} We say that the section $s$ is {\it geometric} if the image $s(G_k)$ of $s$ is contained (hence equal to) in the 
decomposition group $D_x\subset \Pi_X$ associated to a rational point $x\in X(k)$.
\enddefinition

Recall the tower of finite \'etale covers 
$... \to X_{t+1}\to X_{t}\to ...\to X_1=X$ in 3.4.I,
and the section $s:G_k\to \Pi_{X_t}$
of the projection $\Pi_{X_t}\twoheadrightarrow G_k$ induced by $s$, $\forall t\ge 1$.
Assume that the section $s:G_k\to \Pi_X$ has a cycle class uniformly orthogonal to $\Pic$ mod-$p^2$ (cf. Definition 3.4.1).
In particular, the induced section $s:G_k\to \Pi_{X_t}$ also has a cycle class uniformly orthogonal to $\Pic$ mod-$p^2$ (cf. loc. cit.).
There exists, $\forall t\ge 1$, a section 
$$s_t':G_k\to G_{X_t}^{(1/p^2-\sol)}$$ 
of the projection  $G_{X_t}^{(1/p^2-\sol)}  \twoheadrightarrow G_k$ which lifts the section $s$ (cf. Theorem 3.4.4). 
(Note that $k$ satisfies the condition $\bold {(H_{p^2})}$.)
Given integers $t_1\ge t_2\ge1$, we have a commutative diagram
$$
\CD
G_{X_{t_1}}^{(1/p^2-\sol)} @>>> G_k\\
@VVV @V{\id}VV\\
G_{X_{t_2}}^{(1/p^2-\sol)} @>>> G_k\\
\endCD
$$ 
where the left vertical map is induced by the scheme morphism $X_{t_1}\to X_{t_2}$.
We say that the above sections $\{s_t'\}_{t\ge 1}$ are {\bf compatible} if $\forall t_1\ge t_2\ge1$ we have a commutative diagram.
$$
\CD
G_k  @>{s_{t_1}'}>> G_{X_{t_1}}^{(1/p^2-\sol)}\\
@V{\id}VV @VVV\\
G_k  @>{s_{t_2}'}>> G_{X_{t_2}}^{(1/p^2-\sol)}\\
\endCD
$$ 

\proclaim{Lemma 4.2} With the above notations, let $s'=s'_1:G_k\to G_{X}^{(1/p^2-\sol)}$ be a section 
of the projection  $G_{X}^{(1/p^2-\sol)}  \twoheadrightarrow G_k$ which lifts the section $s$. Then $s'$
induces naturally compatible sections $s_t':G_k\to G_{X_t}^{(1/p^2-\sol)}$ of the projections  $G_{X_t}^{(1/p^2-\sol)}  \twoheadrightarrow G_k$ 
which lift the section $s$, $\forall t\ge 1$.
\endproclaim

\demo{Proof} Follows from the fact that we have a commutative diagram
$$
\CD
1@>>> \Cal I_{X_t}  @>>> G_{X_t} @>>> \Pi_{X_t} @>>> 1\\
@. @VVV   @VVV  @VVV\\
1@>>> \Cal I_{X}  @>>> G_{X} @>>> \Pi_{X} @>>> 1\\
\endCD
$$
where the right square is cartesian, and $\Cal I_{X_t}=\Cal I_X$. In particular, $\Cal I_{X_t,1/p^2}=\Cal I_{X,1/p^2}$ and $G_{X_t}^{(1/p^2-\sol)}$ 
is the pull back of the group extension $1\to \Cal I_{X,1/p^2}\to G_{X}^{(1/p^2-\sol)}\to \Pi_X\to 1$ via the natural inclusion $\Pi_{X_t}\hookrightarrow  \Pi_X$, $\forall t\ge 1$.
\qed
\enddemo

Next, recall the exact sequence (cf. diagram (3.6), the case $i=1$)
$$1\to \Cal I_X[p,2]\defeq \underset{U} \to{\varprojlim}\lgroup \underset{j\ge 1} \to{\varprojlim}((\Delta _{2,j}')^U)_{1/p}\rgroup \to G_X^{(p,2)}\to G_{X}^{(1/p^{2}-\sol)}\to 1.$$
We have a commutative diagram
$$
\CD
1@>>> \Cal I_{X_t}[p,2]@>>> G_{X_t}^{(p,2)}@>>> G_{X_t}^{(1/p^{2}-\sol)}@>>> 1\\
@. @VVV   @VVV  @VVV\\
1@>>> \Cal I_X[p,2]@>>> G_X^{(p,2)}@>>> G_{X}^{(1/p^{2}-\sol)}@>>> 1\\
\endCD
$$
where the right square is cartesian (as one easily verifies). In particular, $\Cal I_{X_t}[p,2]=\Cal I_{X}[p,2]$.
Let $s':G_k\to G_{X}^{(1/p^2-\sol)}$ be a section 
of the projection  $G_{X}^{(1/p^2-\sol)}  \twoheadrightarrow G_k$ which lifts the section $s$, and 
$\{s_t':G_k\to G_{X_t}^{(1/p^2-\sol)}\}_{t\ge 1}$ the induced compatible sections 
as in Lemma 4.2 which lift the section $s$.
Let $\tilde s:G_k\to G_X^{(p,2)}$ be a section of the projection $G_X^{(p,2)}\twoheadrightarrow G_k$ which lifts the section $s'$. Then $\tilde s$ induces sections 
$\tilde s_t:G_k\to G_{X_t}^{(p,2)}$ of the projections $G_{X_t}^{(p,2)}\twoheadrightarrow G_k$ which lift the section $s_t'$, $\forall t\ge 1$, and which are compatible in the sense 
that $\forall t_1\ge t_2\ge 1$ integers we have a commutative diagram (cf. above diagram whose right square is cartesian).
$$
\CD
G_k  @>{s_{t_1}'}>> G_{X_{t_1}}^{(p,2)}\\
@V{\id}VV @VVV\\
G_k  @>{s_{t_2}'}>> G_{X_{t_2}}^{(p,2)}\\
\endCD
$$

Also, recall the notations and definitions in 3.4.II, the case ${i=1}$,
relative to the sections $s_{U}':G_k\to \Pi_{U}^{(1/p^2-\sol)}$ induced by $s'$, $\forall$ nonempty open subscheme $U\subseteq X$. Thus, the 
$\{\widetilde X_{j}^{U}\}_{j\ge 1}$ 
are defined in this case $\forall U\subseteq X$ open 
(cf. loc. cit.); they form a system of neighbourhoods of the section $s_{U}'$.
Suppose that the group extension  
$1\to  \pi_1(\widetilde X_{j}^{U}\times_k\overline k,\overline \eta_{j})^{1/p} \to \pi_1(\widetilde X_{j}^{U},\eta_{j})^{(1/p)} \to G_k \to 1$
splits, or equivalently that the class of $J_{j}^1[U]$ in $H^1(G_k,J_{j}[U])$ is divisible by $p$, $\forall j\ge 1$, and  
$\forall U\subseteq X$ as above (cf. Lemma 3.4.8). Then the section $s':G_k\to G_{X}^{(1/p^2-\sol)}$ lifts to a section 
$\tilde s:G_k\to G_{X}^{(p,2)}$
of the projection $G_{X}^{(p,2)}\twoheadrightarrow G_k$ (cf. Theorem 3.4.10).
Moreover, $\tilde s$ induces compatible sections 
$\tilde s_t:G_k\to G_{X_t}^{(p,2)}$ of the projections $G_{X_t}^{(p,2)}\twoheadrightarrow G_k$ which lift the section $s_t'$, $\forall t\ge 1$ (cf. above discussion).
In particular, the above sections $\tilde s_t$ induce naturally sections
$$\rho_t:G_k\to G_{X_t}^{(2/p-\sol)}$$
of the projections $G_{X_t}^{(2/p-\sol)}\twoheadrightarrow G_k$
which lift the sections $s$, $\forall t\ge 1$ (cf. Lemma 3.3.1, as well as the diagrams (3.7) and (3.8)).

\definition{Definition 4.3} With the above notations, we say that the section $s$ is {\bf admissible} if the following two conditions hold.

\noindent
{A1)}\ The section $s$ has a {\bf cycle class uniformly orthogonal to $\Pic$ mod-$p^2$}. 

\noindent
{A2)}\  There {\bf exists} a section $s':G_k\to G_{X}^{(1/p^2-\sol)}$ 
of the projection  $G_{X}^{(1/p^2-\sol)}  \twoheadrightarrow G_k$ which {\bf lifts} the section $s$ (this holds if condition A1 is satisfied by Theorem 3.4.4) 
such that the following holds.
The class of $J_{j}^1[U]$ in $H^1(G_k,J_{j}[U])$ is {\bf divisible} by $p$, $\forall U\subseteq X$ nonempty open subscheme, 
and $\forall j\ge 1$. Or, equivalently, the section $s'$ {\bf lifts} to a section 
$\tilde s:G_k\to G_{X}^{(p,2)}$ of the projection $G_{X}^{(p,2)}\twoheadrightarrow G_k$ (cf. Theorem 3.4.10). 
\enddefinition

\proclaim {Lemma 4.4} Let $k'/k$ be a finite extension, $X_{k'}\defeq X\times_kk'$, and $s_{k'}:G_{k'}\to \Pi_{X_{k'}}$
the section of the projection $\Pi_{X_{k'}}\twoheadrightarrow G_{k'}$ which is induced by $s$. Assume $s$ is admissible then $s_{k'}$ is admissible.
\endproclaim

\demo{Proof} First, if $s$ has a cycle class uniformly orthogonal to $\Pic$ mod-$p^2$ then so does $s_{k'}$ (cf. Definition 3.4.1(ii)). 
The second assertion follows from the various definitions.
\qed
\enddemo

The following is our main result in this section; it provides a  
characterisation of sections of (geometrically pro-$\Sigma$, $p\in \Sigma$) 
arithmetic fundamental groups of $p$-adic curves which are geometric.

\proclaim {Theorem 4.5} We use the above notations. The section $s$ is {\bf admissible} (cf. Definition 4.3)
{\bf if and only if} $s$ is {\bf geometric} (cf. Definition 4.1).
\endproclaim

\demo{Proof} 
The if part follows easily from the various Definitions. We prove the only if part.
 Assume that $s$ is admissible, and $k$ contains a primitive p-th root $\zeta_p$ of $1$. The section $s:G_k\to \Pi_{X_t}$ lifts to a section 
$\rho_t:G_k\to G_{X_t}^{(2/p-\sol)}$ of the projection $G_{X_t}^{(2/p-\sol)}\twoheadrightarrow G_k$ (cf. discussion before Definition 4.3), $\forall t\ge 1$.
The section $\rho_t$ induces a section $\tilde \rho_t: (G_k)_{2/p}\to (G_{X_t})_{2/p}$  of the projection
$(G_{X_t})_{2/p}\twoheadrightarrow ({G_k})_{2/p}$, where the $(\ \ )_{2/p}$ of the various profinite groups
are the second quotients of the $\Bbb Z/p\Bbb Z$-derived series (cf. 1.2). 
The section $\tilde \rho_t$ is geometric and arises from a rational point $x_t\in X_t(k)$ by a result of Pop
(cf. [Pop]). (Here one uses the fact that $\zeta_p\in k$.) In particular, $X_t(k)\neq \varnothing$, $\forall t\ge 1$.  
A well-known limit argument shows that $s$ is geometric (cf. [Tamagawa], Proposition 2.1, (iv), see also the details of the proof of Theorem A in [Sa\"\i di1]). 
In case $\zeta_p\notin k$, let $k'\defeq k(\zeta_p)$. The section $s_{k'}$ is admissible (cf. Lemma 4.4), hence is geometric by the above discussion.
One then verifies easily that $s$ is geometric (cf. [Sa\"\i di2], proof of Theorem B). 
\qed
\enddemo

In the course of proving Theorem 4.5 we proved the following (cf. discussion before Definition 4.3, and the proof of Theorem 4.5).
\proclaim {Proposition 4.6} Let $\tilde s:G_k\to G_X^{(p,2)}$ be a section of the projection $G_X^{(p,2)}\twoheadrightarrow G_k$,
and $s:G_k\to \Pi_X$ the section of the projection $\Pi_X\twoheadrightarrow G_k$ which is induced by $\tilde s$. Then $s$ is geometric. 
\endproclaim

\definition {Remarks 4.7} {1)} Theorem 4.5 above is stronger and more precise than Theorem A in [Sa\"\i di1].

\noindent
{2)}\ There are examples of sections $s:G_k\to \Pi_X$ as above, where $\Sigma =\{p\}$, which are {\it not} geometric (cf. [Hoshi]). 
These provide examples of sections $s$ as above which are {\it not} admissible by Theorem 4.5 (where $\Sigma=\{p\}$).
It would be interesting to know which of the conditions A1 and A2 in the definition of admissible sections fail to hold in Hoshi's example.
In [Sa\"\i di3] we observe that the section in Hoshi's example is orthogonal to $\Pic^0$ in the sense that the map
$s^{\star}:H^2(\Pi_X,\Bbb Z_p)\to H^2(G_k,\Bbb Z_p)$ annihilates the image of $\Pic^{0}(X)$.
\enddefinition

The following is an application of our results to the absolute anabelian geometry of $p$-adic curves.
\proclaim{Theorem 4.8} Let $p_X,p_Y\in \Primes$, and
$X$ (resp. $Y$) a proper smooth and geometrically connected hyperbolic curve over a $p_X$-adic local field
$k_X$ (respectively,  $p_Y$-adic local field $k_Y$). Let $p_X\in \Sigma_X$ (resp. $p_Y\in \Sigma_Y$) be a non-empty set of prime integers of cardinality $\ge 2$, 
$\Pi_X$ (resp. $\Pi_Y$) the geometrically pro-$\Sigma_X$ (resp. pro-$\Sigma_Y$) arithmetic fundamental group of $X$ (resp. $Y$),
and $\varphi:\Pi_X\to \Pi_Y$ an isomorphism of profinite groups which fits in the following commutative diagram
$$
\CD
G_X^{(p_X,2)} @>{\widetilde \varphi}>> G_Y^{(p_Y,2)}\\
@VVV   @VVV\\
\Pi_X@>{\varphi}>>\Pi_Y\\
\endCD
$$
where $\widetilde \varphi$ is an isomorphism of profinite groups, and the vertical maps are the natural projections. Then $\varphi$ is geometric, i.e., arises from a uniquely determined isomorphism of schemes $X\isom Y$.
\endproclaim

\demo{Proof}
The existence of the lifting $\widetilde \varphi$ of $\varphi$ implies, by Proposition 4.6, that $\varphi$ preserves the decomposition groups at closed points. The statement follows then from [Mochizuki1], Corollary 2.9.
\qed
\enddemo

\subhead
\S 5. Local sections of arithmetic fundamental groups of $p$-adic curves 
\endsubhead
We prove that a certain class of sections of arithmetic fundamental groups of $p$-adic curves are (uniformly) orthogonal to $\Pic^{\wedge}$. 
We use the notations in $\S4$.
\subhead
5.1. Arithmetic fundamental groups of formal fibres of $p$-adic curves
\endsubhead
Let $\Cal O_k$ be the valuation ring of $k$, and $\widetilde X\to \Spec \Cal O_k$ a flat and proper model of $X$ over $\Cal O_k$ with $\widetilde X$ normal. Let $x\in \widetilde  X^{\cl}$ be a closed point, and $\hat \Cal O_{\widetilde X,x}$ the completion of the local ring of $\widetilde X$ at $x$. We will refer to $\Cal X\defeq \Spec (\hat \Cal O_{\widetilde X,x}\otimes _{\Cal O_k}k)$ as the {\it formal fibre} of $\widetilde X$ at $x$ (or simply a formal fibre). Assume $\Cal X$ is geometrically connected, and write $\overline {\Cal X}\defeq \Cal X\times _k\overline k$.
Let $\bar \beta$ be a geometric point of $\overline {\Cal X}$, which determines a geometric point $\beta $ of ${\Cal X}$. 
Write $\Delta_{\Cal X}\defeq \pi_1(\overline {\Cal X},\overline \beta)^{\Sigma}$
for the maximal pro-$\Sigma$ quotient of $\pi_1(\overline {\Cal X},\overline \beta)$,
and
$\Pi_{\Cal X}\defeq  \pi_1(\Cal X, \beta)/ \Ker  (\pi_1(\overline {\Cal X},\overline \beta)\twoheadrightarrow
\pi_1(\overline {\Cal X},\overline \beta)^{\Sigma})$. 
We have a commutative diagram of exact sequences 
$$
\CD
1 @>>> \Delta_{\Cal X}@>>> \Pi_\Cal X @>>> G_k @>>> 1\\
@.  @VVV  @VVV  @V{\id}VV\\
1 @>>> \Delta_X@>>> \Pi _X@>>> G_k @>>> 1\\
\endCD
\tag 5.1$$
where the middle vertical map (defined up to inner conjugation)
is induced by the scheme morphism $\Cal X\to X$.

For the rest of this section we assume that $\Cal X$ is a formal fibre as in 5.1, which is geometrically connected.
\definition {Definition 5.2} A section $\tilde s:G_k\to \Pi_{\Cal X}$ of the projection $\Pi_{\Cal X}\twoheadrightarrow G_k$ induces a section $s:G_k\to \Pi_X$ of the projection 
$\Pi_X\twoheadrightarrow G_k$ (cf. diagram (5.1)). We will refer to such a section $s$ as a {\bf local section} of the projection $\Pi_X\twoheadrightarrow G_k$.
\enddefinition

Note that a geometric section (cf. Definition 4.1) is a local section in the above sense, as one easily verifies. Our main result is the following.

\proclaim {Theorem 5.3} Let $s:G_k\to \Pi_X$ be a {\bf local} section of the projection $\Pi_X\twoheadrightarrow G_k$. Then $s$ has {\bf a cycle class uniformly orthogonal to $\Pic^{\wedge}$} in the sense of [Sa\"\i di], Definition 1.4.1(i). 
\endproclaim

\demo{Proof} Let $\Cal X=\Spec (\hat \Cal O_{\widetilde X,x}\otimes _{\Cal O_k}k)$ be as in 5.1, and $\tilde s:G_k\to \Pi_{\Cal X}$ a section of the projection $\Pi_{\Cal X}\twoheadrightarrow G_k$ which induces the section $s:G_k\to \Pi_X$. 
Let $\{X_t\}_{t\ge 1}$ be as in $\S4$, $s_t:G_k\to \Pi_{X_t}$ the section of the projection $\Pi_{X_t}\twoheadrightarrow G_k$ which is induced by $s$, and $s_t^{\star}:H^2(X_t,\hat \Bbb Z(1)^{\Sigma})
\to H^2(G_k,\hat \Bbb Z(1)^{\Sigma})$ the retraction map induced by $s_t$, $\forall t\ge 1$. We show $s_t^{\star}(\Pic(X_t)^{\wedge})=0$ (cf. loc. cit. for the definition of $\Pic(X_t)^{\wedge}$). 
Recall the continuous homomorphism 
$\phi:\Pi_{\Cal X}\to \Pi_X$ (cf. the middle vertical map in diagram (5.1)).
Then $\Pi_{\Cal X_t}\defeq \phi^{-1}(\Pi_{X_t})$ is an open subgroup containing $\tilde s(G_k)$, and corresponds to an \'etale cover $\Cal X_t\to \Cal X$ with $\Cal X_t$ geometrically connected (as $\Pi_{\Cal X_t}$
projects onto $G_k$ via the projection $\Pi_{\Cal X}\twoheadrightarrow G_k$). Moreover, the section $\tilde s$ induces a retraction $\tilde s_t^{\star}:H^2(\Pi _{\Cal X_t},\hat \Bbb Z(1)^{\Sigma})
\to H^2(G_k,\hat \Bbb Z(1)^{\Sigma})$ of the natural map $H^2(G_k,\hat \Bbb Z(1)^{\Sigma})\to H^2(\Pi_{\Cal X_t},\hat \Bbb Z(1)^{\Sigma})$ induced by the projection 
$\Pi_{\Cal X_t}\twoheadrightarrow G_k$. Let $A\defeq \hat \Cal O_{\widetilde X,x}$,
and $A_k\defeq \hat \Cal O_{\widetilde X,x}\otimes _{\Cal O_k}k$. The open subgroup $\Pi_{\Cal X_t}$ corresponds to an \'etale cover $\Cal X_t=\Spec B_k\to \Cal X=\Spec A_k$,
where $B_k/A_k$ is an \'etale extension. Let $B$ be the integral closure of $A$ in $B_k$. Thus, $B$ is a complete local ring of dimension 2, which dominates $\Cal O_k$, and the residue field of $B$ is finite. We have a scheme theoretic morphism $\Cal X_t\to X_t$. 
Further, we have an injective homomorphism
$H^2(\Pi_{\Cal X_t},\hat \Bbb Z(1)^{\Sigma})\hookrightarrow H_{\et}^2(\Cal X_t,\hat \Bbb Z(1)^{\Sigma})$ arising from the Cartan-Leray spectral sequence
(cf. [Serre], proof of Proposition 1), as well as an injective Kummer homomorphism $\Pic (\Cal X_t)^{\wedge}\hookrightarrow 
H_{\et}^2(\Cal X_t,\hat \Bbb Z(1)^{\Sigma})$, where $\Pic^{\wedge}\defeq \Pic\otimes_{\Bbb Z} \hat \Bbb Z^{\Sigma}$. On the other hand we have a commutative diagram of homomorphisms
$$
\CD
H^2(\Pi_{\Cal X_t},\hat \Bbb Z(1)^{\Sigma}) @>{\tilde s_t^{\star}}>>  H^2(G_k,\hat \Bbb Z(1)^{\Sigma}) \\
@A\psi_tAA  @A{\id}AA\\
H^2(\Pi_{X_t},\hat \Bbb Z(1)^{\Sigma}) @>{s_t^{\star}}>>  H^2(G_k,\hat \Bbb Z(1)^{\Sigma}) \\
\endCD
$$
where $\psi_t$ is induced by the map $\phi:\Pi_{\Cal X_t}\to \Pi_{X_t}$.

We claim that $\psi_t(\Pic(X_t)^{\wedge})$ is {\it torsion}, from this it follows that $s_t^{\star}(\Pic(X_t)^{\wedge})=\tilde s_t^{\star}(\psi_t(\Pic(X_t)^{\wedge}))=0$, since
$H^2(G_k,\hat \Bbb Z(1)^{\Sigma})\isom \hat \Bbb Z^{\Sigma}$ is torsion free. Indeed, we have a pull-back morphism $\Pic (X_t)^{\wedge}\to \Pic(\Cal X_t)^{\wedge}$, 
which fits in the commutative diagram
$$
\CD
\Pic(\Cal X_t)^{\wedge} @>>>  H_{\et}^2(\Cal X_t,\hat \Bbb Z(1)^{\Sigma}) \\
@A{\id}AA   @AAA \\
\Pic(\Cal X_t)^{\wedge} @.  H^2(\Pi_{\Cal X_t},\hat \Bbb Z(1)^{\Sigma}) \\
@AAA   @A{\psi_t}AA \\
\Pic(X_t)^{\wedge} @>>>  H^2(\Pi_{X_t},\hat \Bbb Z(1)^{\Sigma}) \\
\endCD
$$
where the horizontal maps are injective Kummer maps (recall the identification  $H^2(\Pi_{X_t},\hat \Bbb Z(1)^{\Sigma})\isom H_{\et}^2(X_t,\hat \Bbb Z(1)^{\Sigma})$), and the upper right vertical map is the injective map discussed above. Our claim follows then from the following.

\proclaim {Proposition 5.4} With the notations above $\Pic(\Cal X_t)$, and a fortiori $\Pic(\Cal X_t)^{\wedge}$, is {\bf finite}.
\endproclaim

\demo{Proof of Proposition 5.4}
This follows from the fact, proven by Shuji Saito, that $\Pic(\Spec  B\setminus \{m_B\})$ is finite, where $B$ 
is as in the above discussion and $m_B$ is its maximal ideal (cf. [Saito], Theorem 0.11). 
\qed
\enddemo
This finishes the proof of Theorem 5.3.
\qed
\enddemo

Finally, we provide the following characterisation of local sections which are geometric. We use the above notations.

\proclaim{Theorem 5.5} Let $s:G_k\to \Pi_X$ be a {\bf local} section of the projection $\Pi_X\twoheadrightarrow G_k$ (cf. Definition 5.2). Then $s$ {\bf lifts} to a section $\rho:G_k\to G_X^{(1-\sol)}$
(resp. $\rho_n:G_k\to G_X^{(1/p^n-\sol)}$) of the projection $G_X^{(1-\sol)}\twoheadrightarrow G_k$ (resp. $G_X^{(1/p^n-\sol)}\twoheadrightarrow G_k$, $\forall n\ge 1$). Moreover,
the section $s$ is {\bf geometric if and only if} there exists a lifting of $s$ to a section $\rho_2:G_k\to G_X^{(1/p^2-\sol)}$ as above such that
one of the following equivalent conditions hold.

\noindent
{(i)} With the notations in $\S4$ (cf. the discussion before Definition 4.3), 
the class of $J_{j}^1[U]$ in $H^1(G_k,J_{j}[U])$ is {\bf divisible} by $p$, $\forall U\subseteq X$ nonempty open subscheme, 
and $\forall j\ge 1$.

\noindent
{(ii)} The section $\rho_2$ {\bf lifts} to a section 
$\tilde \rho_2:G_k\to G_{X}^{(p,2)}$ of the projection $G_{X}^{(p,2)}\twoheadrightarrow G_k$. 

\endproclaim

\demo{Proof} The first assertion follows from Theorem 5.3, Theorem 3.4.4, and Theorem 2.3.5 in [Sa\"\i di]. The second assertion follows from Theorem 3.4.10, and Theorem 4.5.
\qed
\enddemo

$$\text{References}$$

\noindent
[Grothendieck] Grothendieck, A., Rev\^etements \'etales et groupe fondamental, Lecture 
Notes in Math. 224, Springer, Heidelberg, 1971.

\noindent
[Harari-Szamuley], Harari, D., Szamuely, T., Galois sections for abelianized fundamental groups, Math. Ann. (2009) 344, 779-800.

\noindent
[Hoshi] Hoshi, Y., Existence of nongeometric pro-$p$ Galois sections of hyperbolic curves, Publ. Res. Inst. Math. Sci. 46 (2010), no. 4, 829-848.

\noindent
[Mochizuki] Mochizuki, S., Absolute anabelian cuspidalizations of proper hyperbolic curves,  J. Math. Kyoto
Univ.  47  (2007),  no. 3, 451--539.

\noindent
[Mochizuki1] Mochizuki, S., Topics in absolute anabelian geometry II: decomposition groups and endomorphisms.
J. Math. Sci. Univ. Tokyo, Vol. 20 (2013), No. 2, Page 171-269.

\noindent
[Pop] Pop, F., On the birational $p$-adic section Conjecture, Compos. Math. 146 (2010), no. 3, 621-637.

\noindent
[Sa\"\i di] Sa\"\i di, M., The cuspidalisation of sections of arithmetic fundamental groups, Advances in Mathematics, 230 (2012) 1931-1954.

\noindent
[Sa\"\i di1] Sa\"\i di, M., On the $p$-adic section conjecture, Journal of Pure and Applied Algebra, Volume 217, Issue 3, (2013), 583-584.

\noindent
[Sa\"\i di2] Sa\"\i di, M., On the section conjecture over function fields and finitely generated fields, Publ. RIMS Kyoto Univ. 52 (2016), 335-357.

\noindent
[Sa\"\i di3] Sa\"\i di, M., Arithmetic of $p$-adic curves and sections of geometrically abelian fundamental groups. Manuscript.


\noindent
[Saito] Saito, S., Arithmetic on two dimensional local rings, Invent. math. 85, 379-414 (1986).

\noindent
[Serre] Serre J.-P., Construction de rev\^etements de la droite affine en caract\'eristique $p>0$, C. R. Acad. Sci. Paris 311 (1990), 341-346.

\noindent
[Tamagawa] Tamagawa, A., The Grothendieck conjecture for affine curves,  Compositio Math.  109  (1997),  no. 2, 135--194.

\bigskip

\noindent
Mohamed Sa\"\i di

\noindent
College of Engineering, Mathematics, and Physical Sciences

\noindent
University of Exeter

\noindent
Harrison Building

\noindent
North Park Road

\noindent
EXETER EX4 4QF 

\noindent
United Kingdom

\noindent
M.Saidi\@exeter.ac.uk

\end
\enddocument